\documentclass[a4paper,11pt]{article}

\usepackage[T1]{fontenc}
\usepackage{amsmath}
\usepackage{amssymb}
\usepackage{amsthm}
\usepackage{amscd}
\usepackage{url}
\usepackage{euscript}
\usepackage{mathrsfs}
\usepackage{relsize} 
\usepackage{dsfont}
\usepackage{color}

\usepackage{hyperref} 

\usepackage{bbm}

\numberwithin{equation}{section}  

\newtheorem{theorem}{{\scshape {\bfseries Theorem}}}[section]
\newtheorem{lemma}[theorem]{{\scshape {\bfseries Lemma}}}
\newtheorem{proposition}[theorem]{{\scshape {\bfseries Proposition}}}
\newtheorem{corollary}[theorem]{{\scshape {\bfseries Corollary}}}
\newtheorem{remark}{{\scshape {\bfseries Remark}}}[section]

\renewenvironment{proof}[1]
{\noindent{{\bf{\small{P}{\scriptsize ROOF}}}.}\hspace{0.1cm} #1}
{$\;\blacksquare$\newline}

\newcommand{\entiers}{\mathds{N}^{*}}

\newcommand{\Oun}{\mathcal{O}(1)} 
\newcommand{\geno}{\tilde L}
\newcommand{\gen}{{L}}
\newcommand{\un}{\mathbbm{1}}
\newcommand{\proba}{{\mathbb P}}
\newcommand{\real}{\mathds{R}}
\newcommand{\integers}{\mathds{N}}
\newcommand{\dir}{\mathscr{E}}
\newcommand{\esperance}{{\mathds{E}}}
\newcommand{\ns}{n_{*}{\scriptstyle (K)}}
\newcommand{\nss}{n_{**}{\scriptscriptstyle (K)}}
\newcommand{\nsss}{n_{***}{\scriptscriptstyle (K)}}

\newcommand{\tpg}[1]{The proof of this theorem is given in Section \ref{#1}}

\font\gfont=cmmi10 scaled \magstep{1.5}     
\newcommand{\gdelta}{\hbox{\gfont \char14}}
\newcommand{\grho}{\hbox{\gfont \char26}}

\title{Sharp asymptotics for the quasi-stationary distribution of birth-and-death processes}

\author{J.-R. Chazottes$^{\textup{(a)}}$, P. Collet$^{\textup{(a)}}$, S. M\'el\'eard$^{\textup{(b)}}$\\
{\small $^{\textup{(a)}}$ Centre de Physique Th\'eorique, CNRS UMR 7644, Ecole polytechnique}\\
{\small F-91128 Palaiseau Cedex (France)}\\
{\small $^{\textup{(b)}}$ Centre de Math\'ematiques Appliqu\'ees, CNRS UMR 7641, Ecole polytechnique}\\
\small{F-91128 Palaiseau Cedex (France)}
}

\begin{document}

\maketitle

\begin{abstract}
We study a general class of birth-and-death processes with state space $\integers$ that describes the size of a
population  going to extinction with probability one. This class contains the logistic case.
The scale of the population is measured in terms of a `carrying capacity' $K$.
When $K$ is large, the process is expected to stay close to its deterministic equilibrium during a long time
but ultimately goes extinct. 
Our aim is to quantify the behavior of the process and the mean time to extinction in the quasi-stationary distribution as a function of $K$, for large $K$. We also give a quantitative description of this quasi-stationary distribution. It turns out to be close to a Gaussian distribution centered about the deterministic long-time equilibrium, when $K$ is large. \\
Our analysis relies on precise estimates of the maximal eigenvalue, of the corresponding eigenvector and of the spectral gap of a self-adjoint operator associated with the semigroup of the process.
\end{abstract}

\newpage

\tableofcontents

\newpage

\section{Introduction}

We study a general class of birth-and-death processes with state space $\integers$ that describes the size of a
population going to extinction with probability one. 
For a population of size $n\in \integers^*$, the  birth rate is denoted by $\lambda_n>0\,$ and the death
rate by $\mu_n>0$.  Furthermore, we assume that
\[
\lambda^{\scriptscriptstyle{K}}_{n}=n\;\tilde \lambda\left(\frac{n}{K}\right),\quad
\mu^{\scriptscriptstyle{K}}_{n}=n\;\tilde \mu\left(\frac{n}{K}\right)
\]
where $\tilde \lambda,\tilde \mu$ are positive functions and $K$ is a scaling parameter describing the amount of available resources (that is called  the `carrying capacity' in ecology). We assume that $\lambda_0=\mu_0=0$, entailing absorption at state $0$. 

In this work, we consider the case where absorption at $0$ happens with probability one.
We also assume that the time to this absorption has finite expectation. In this situation, the unique stationary probability measure is $\gdelta_0$, the Dirac mass at state $0$. In order to understand the behavior of the process before absorption,  a relevant object to look at is a so-called quasi-stationary distribution, \emph{i.e}, a probability distribution that is stationary when the process is conditioned to survive. Our aim is to describe what happens for large $K$.

The prominent example is the so-called logistic birth-and-death process
$(X_t^{\scriptscriptstyle{K}},t\geq 0)$ defined by following birth and death rates
\begin{equation}\label{logistic-ex}
\lambda^{\scriptscriptstyle{K}}_n=\tilde{\lambda}\, n,\; 
\mu^{\scriptscriptstyle{K}}_n=n\left(\tilde{\mu}+\frac{n}{K}\right)
\end{equation}
for $n\geq 1$, where $\tilde \lambda, \tilde \mu$ are positive parameters. 
It is a classical result (see {\em e.g.} \cite{K1}) that if the process starts in a state of the form $\lfloor x_0 K\rfloor$ ($x_0>0$), then the rescaled process $X_t^{\scriptscriptstyle{K}}/K$ is `close',  in the limit  as $K\to\infty$, during any given finite interval of time, to the solution of the differential equation
\[
\frac{\mathrm{d}x}{\mathrm{d}t}=x(\tilde{\lambda}-\tilde{\mu}-{x})
\]
with initial condition $x_0$. 
This differential equation has a unique attractive equilibrium $x_*=\tilde{\lambda}-\tilde{\mu}$ and  the integer $\lfloor x_* K\rfloor$ can be considered as an approximation of the population size over every given finite time interval.
However, for each $K$, the process $X_t^{\scriptscriptstyle{K}}$ goes almost surely to extinction as $t\to\infty$, see \cite{VD}. 

In this paper, we consider more general processes with the same kind of behavior. One of our motivations is
to quantify, as a function of $K$, the scale of the mean time to extinction, the time-scale of convergence to the quasi-stationary distribution, and the time-scale during which the process is close to the rescaled deterministic equilibirum
$\lfloor x_* K\rfloor$ with high probability.  

Our results can be colloquially described as follows. 
We get an upper bound of order $K\log K$ for the time it takes for the process to be close to the quasi-stationary distribution. 
We also get the existence of a time interval, exponentially long in $K$, during which the process, if we start from a population of order $K$, is nearly  distributed according to the quasi-stationary distribution.\\
We also prove that the total variation distance between the quasi-stationary distribution and a  Gaussian distribution is bounded by $1/\sqrt{K}$. This Gaussian distribution is centered around $\lfloor x_* K\rfloor$ and its variance is of order $K$.\\
As a by-product of our analysis we show that the mean time to extinction with respect to the quasi-stationary distribution is given by 
\[
\frac{c}{\mathsmaller{\sqrt{K}}}\;\mathlarger{e}^{K\mathlarger{\int}_0^{x_*} \log \frac{\tilde{\lambda}(x)}{\tilde{\mu}(x)}  \,\mathrm{d} x}\,
\left( 1+ \mathsmaller{\mathcal{O}\left(\frac{(\log K)^3}{\sqrt{K}}\right)}\right)
\]
where $c$ is a constant independent of $K$ that is explicitly given later on. Roughly speaking,
this mean time is exponentially large in $K$.

Motivated by population extinction in biology, many people attempted to analyze quasi-stationary distributions. But even in the simplest models, like the logistic model, this turned out to be a complicated task. Previous results are mostly based on either Monte-Carlo simulations or uncontrolled approximations based on heuristic ansatzes, see the review paper \cite{OM} and also \cite{nasell,KS}. The present work is the first one in which controlled mathematical approximations are obtained for the quasi-stationary distribution for a class of models encompassing the logistic model. 

We are aware of only a few mathematical results related to our work. In \cite{DSS}, the authors do not study the quasi-stationary distribution but only the mean time to expectation starting from a state of order $K$ for which they obtain the asymptotic behavior in $K$ (see also \cite{SS}).
Here we are able to control this quantity for all initial states and also for the quasi-stationary distribution as a starting distribution.
In \cite{barbour}, the authors show that the quasi-stationary distribution can be approximated in total variation distance by an auxiliary process called the `returned process'. They also prove a bound for the total variation distance between the law of the process $X_t^{\scriptscriptstyle{K}}$ for fixed values of $t$ and the quasi-stationary distribution. This is somewhat related to one of our theorems (Theorem \ref{thm:dtv}). 
Let us also mention the articles \cite{CV,CT,DM} about quantitative convergence to quasi-stationarity.

The main tool in this work is the analysis of an operator $L$ that is related to the generator of the killed process. We use a weighted Hilbert space  where $L$ is self-adjoint. The operator $L$ has a maximal simple and negative  eigenvalue $-\grho_0$. The mean time to extinction is exactly $1/\grho_{0}$. The quasi-stationary distribution is constructed  from the corresponding positive eigenvector.The method of  analysis of the equation $Lu=-\grho_0 u$ is inspired by matching techniques reminiscent of the WKB method in Physics
\cite{fedor,levinson}.

\section{Standing assumptions and notations}\label{sec:hyp}

In the sequel most quantities will depend on the parameter $K$. We will not indicate systematically this dependence in the notation,
except when we want to highlight it. 
Recall that 
\begin{equation}\label{eq:scaling}
\lambda_{n}=n\;\tilde \lambda\left(\frac{n}{K}\right),\;
\mu_{n}=n\;\tilde \mu\left(\frac{n}{K}\right).
\end{equation}

In the rest of the paper, the functions $x\mapsto\tilde{\lambda}(x)$ and $x\mapsto\tilde{\mu}(x)$, defined on $\real_+$, are assumed to be positive, differentiable and increasing. 
In particular, this implies that the sequences $(\lambda_n)_n$ and $(\mu_n)_n$ are increasing. 

From now on, we assume that the following properties for the functions $\tilde{\lambda}$ and $\tilde{\mu}$ hold
throughout the paper.
\begin{align}
\label{infini}
&\bullet\;  \lim_{x\to +\infty} \frac{\tilde \lambda(x)}{\tilde \mu(x)}= 0;\\
&\bullet\;
\label{coeff}
\tilde \lambda(0) > \tilde \mu(0) >0;\\
&\bullet\; 
\label{eq:xstar}
\text{There exists a unique}\;x_{*}\in \real_{+}\;\text{such that}\;\tilde{\lambda}(x_*)=\tilde{\mu}(x_*);\\
&\bullet\;
\label{eq:generic}
\tilde{\lambda}'(x_*)\neq \tilde{\mu}'(x_*);\\
&\bullet\;
\label{cv-int}
\int_{\frac{x_*}{2}}^{+\infty} \frac{\textup{d}x }{x\, \tilde \mu(x)}<+\infty;\\
&\bullet\;
\label{pierre}
\sup_{x\in\real_+} \frac{\tilde{\mu}'(x)}{\tilde{\mu}(x)} <+\infty;\\
&\bullet\;
\nonumber
\text{The function}\;x\mapsto \log \frac{\tilde{\mu}(x)}{\tilde{\lambda}(x)}\;\text{defined on}\;\real_+\;\text{is increasing.}\\
& 
\nonumber
\text{The function}\; H:\real_+\to\real\;\text{defined by}\\
&
\label{def:H}
\qquad\qquad \qquad \qquad H(x)=\int_{x_*}^x \log\frac{\tilde{\mu}(s)}{ \tilde{\lambda}(s)}\, \mathrm{d}s\\
& 
\nonumber
\text{is assumed to have the following properties:}\\
\nonumber
& \bullet \;H\;\text{is}\;\text{three times differentiable;}\\ 
&\bullet\; \sup_{x\in\real_+} (1+x^2) |H'''(x)|< +\infty.
\label{seconde}
\end{align}

Some comments are in order about the above assumptions.
The relevant assumptions from a biological viewpoint are assumptions \eqref{infini}, \eqref{coeff} and \eqref{eq:xstar}. The first one means that, when the population size gets large, deaths prevail. The second one means the opposite: at low population size, births prevail. The third one means that there is a unique equilibrium for the associated differential equation. This rules out for instance the so-called Allee effect where there are two non-trivial equilibria. Assumption \eqref{eq:generic} is a genericity property. The remaining assumptions are technical but they are by far true in the logistic case and in many other models.

We shall denote by $(X_t^{\scriptscriptstyle{K}}, t\geq 0)$ the birth-and-death process associated with the
rates $(\lambda_n)$ and $(\mu_n)$.
Thorough the paper we will use the classical notation 
\begin{equation}\label{eq:pin}
\pi_{n}=\frac{\lambda_1\cdots\lambda_{n-1}}{\mu_1\cdots \mu_n},\;\textup{for}\; n\geq 2
\end{equation}
and we set $\pi_1:=\frac{1}{\mu_1}$. The following trivial identity will be used repeatedly.
\begin{equation}
\label{tricky}
\lambda_{n}\, \pi_{n}= \mu_{n+1}\,\pi_{n+1}.
\end{equation}
  
One can verify that condition \eqref{infini}, together with the facts that $(\mu_n)_n$ is increasing and that
$\tilde{\mu}(0)$ is bounded away from zero, imply the following two properties:
\begin{equation}
\label{extinction}
(\star)\;\mathlarger{\sum}_{n\geq 1} \frac{1}{\lambda_{n}\ \pi_{n}}= +\infty
\quad\text{and}\quad
(\star\star)\; \sum_{n\geq 1} \pi_{n}<+\infty.
\end{equation} 
The property $(\star)$ implies absorption of the process at state $0$ with probability one.
The property $(\star\star)$ ensures finiteness of the expectation of the absorption time, that is,
$\esperance_m[T_{0}]<+\infty$ for every $m\in\entiers$, where $T_0=\inf\{t\geq 0 : X_t^{\scriptscriptstyle{K}}=0\}$.  
We refer to \cite[p. 384]{KT} and \cite[chapter 3]{Allen2011} for details.

Condition \eqref{cv-int} implies  
\begin{equation}
\label{QSD}
\mathlarger{\sum}_{n\geq 1} \frac{1}{\lambda_{n}\ \pi_{n}}
\left(\mathlarger{\sum}_{i\geq n+1} \pi_{i}\right) <+\infty.
\end{equation}
(See Lemma \ref{controle} for a proof.)
As proved in \cite{VD}, this is a sufficient condition for the existence and uniqueness of a quasi-stationary distribution.
It turns out that it is a necessary condition as well as it can be deduced from \cite{cattiaux}.
Condition \eqref{pierre} implies 
\begin{equation}\label{eq:mumu}
\sup_{n} \frac{\mu_{n+1}}{\mu_n}<\infty.
\end{equation}
This follows from the mean value theorem to the function $x\mapsto \log\tilde{\mu}(x)$.
We will assume that 
\begin{equation}\label{eq:pimu}
\pi_n \mu_n^2\xrightarrow{n\to\infty} 0.
\end{equation}
This is a technical condition that we use in the spectral theory of the operator associated with
the process. 

Finally, let us recall (see {\em e.g.} \cite{K1}) that for large $K$,  the process
$(X_t^{\scriptscriptstyle{K}}/K, t\geq 0)$  is close to the solution of the ordinary differential equation
\begin{equation}\label{eq:ODE}
\frac{\mathrm{d}x}{\mathrm{d}t}=x\, \big(\tilde{\lambda}(x)-\tilde{\mu}(x)\big)
\end{equation}
during any given finite time interval.
Our assumptions imply that the differential equation \eqref{eq:ODE} has the unique non-zero equilibrium $x_{*}$.
Observe that, because of assumptions on the functions $x\mapsto \tilde \lambda(x)$ and $x\mapsto \tilde \mu(x)$, one has $\frac{\tilde \lambda(x)}{\tilde \mu(x)} >1$ for $x<x_{*}$ and $\frac{\tilde \lambda(x)}{\tilde \mu(x)}<1$ for $x>x_{*}$. This implies the stability of the equilibrium $x_*$ of the deterministic equation \eqref{eq:ODE} and, using \eqref{eq:generic}, we get
\begin{equation}\label{Hseconde}
H''(x_*)>0.
\end{equation}
We shall use the notation
\begin{equation}\label{eq:nstar}
\ns =\lfloor x_* K\rfloor.
\end{equation}
This quantity plays a natural role in the sequel.

\textbf{An example}. For the logistic birth-and-death process defined in \eqref{logistic-ex}, we have
$\tilde{\lambda}(x)=\tilde{\lambda}$ and $\tilde{\mu}(x)=\tilde{\mu}+x$. If $\tilde{\lambda}>\tilde{\mu}$, 
it is easy to check that all the above conditions are fullfilled.
One has $\ns=\lfloor(\tilde\lambda-\tilde\mu) K\rfloor$. 

\section{Statements of the main results}

\subsection{The generator and its spectrum}

Our goal is to link the semigroup of the process $(X_t^{\scriptscriptstyle{K}},t\geq 0)$ `killed' at $0$ to a self-adjoint operator with compact resolvent in an appropriate Hilbert space. The spectral theory for this operator lies at the core of our work.

Let us denote by $\mathscr{D}$ the set of sequences with finite support on $\entiers$. Define the operator 
$\geno$ with domain $\mathscr{D}$ by
\[
(\geno u)_{n}=\lambda_{n}u_{n+1}+\mu_{n}u_{n-1}\un_{\{n\ge2\}}-(\lambda_{n}+\mu_{n})u_{n}.
\]
We introduce the following weighted space of sequences of complex numbers
\[
\ell^2(\pi) = \Big\{u=(u_n)_{n\in\entiers} :\sum_{n=1}^\infty \pi_{n} |u_{n}|^2 < \infty\Big\}
\]
where the $\pi_n$'s are defined in \eqref{eq:pin}.
The space $\ell^2(\pi)$ is a Hilbert space when endowed with the scalar product
\[
\langle u, v\rangle_{\pi}= \sum_{n=1}^\infty \pi_{n} \bar u_{n} v_{n}
\]
where $ \bar u_{n}$ is the complex conjugate of $u_n$.
We shall denote by $\|\!\cdot\!\|_{\pi}$ the associated norm.

The main content of the following theorem is that one can extend the operator $\tilde \gen$ to an operator $L$
that is the infinitesimal generator of a positive and contractive semigroup in $\ell^2(\pi)$. Moreover, this operator has a discrete spectrum with a maximal eigenvalue that is simple and negative.

\begin{theorem}[The operator $L$, $\grho_0$, $\varphi$ and $\grho_1$]
\label{autoadj}
\leavevmode
\begin{enumerate}
\item 
The operator $\tilde \gen$ is symmetric on $\mathscr{D}$. It  is closable in $\ell^{2}(\pi)$. 
\item
We will denote by $\gen$ its closure and by $\mathcal{D}$ the domain of this closure.
The operator $\gen$ defines a positive  contraction semigroup in $\ell^{2}(\pi)$.
\item
$\gen$ is a dissipative, self-adjoint operator with a compact resolvent.  
Its spectrum is discrete  and  the maximal  eigenvalue is  simple and negative.
We denote it by $-\grho_{0}$. The corresponding eigenvector can be chosen positive
and we denote it by $\varphi$. Finally, we denote by $-\grho_{1}$ the second largest eigenvalue.
\end{enumerate}
\end{theorem}

\tpg{proof-autoadj}.

\begin{remark}
The construction of $\mathcal{D}$ is general; see \cite[III.5.3]{Kato}.
\end{remark}

For all $t>0$, $n,m\in \entiers$, let
\begin{equation}\label{def:Pt}
P_{t}(n,m)= \frac{1}{\pi_{n}} \langle \mathrm{e}_{n}, \mathlarger{e}^{t\gen}  \mathrm{e}_{m}\rangle_{\pi},
\end{equation}
where for each $n$, $\mathrm{e}_{n}$ is defined by $\mathrm{e}_{n}(k) = \delta_{n,k}$ for $k=1,2,\ldots$.
A straightforward computation shows that the `matrix' $(P_{t}(m,n))_{(m,n)\in \entiers\! \times\entiers}$
is a solution of the Kolmogorov equation
\[
\frac{\textup{d}P_t(n,m)}{\textup{d}t}   = \big(\gen P_{t}\big)(n,m) = \big(P_{t}L\big)(n,m).
\]
Furthermore, one can verify that there exists some $M\geq 1$ such that for all $t$ and all $n$,
$\left| \sum_{k=1}^\infty  P_t(n,k) \right| \leq M.$ 
The uniqueness of such a family has been proven in \cite[Theorem 14 p. 528]{KM}  under Assumption \eqref{extinction}.
This implies that the symmetric sub-markovian semigroup $(P_{t},t\geq 0)$ is the extension of the transition semigroup of the Markov process $(X^{\scriptscriptstyle{K}}_{t}, t\geq 0)$ to $\ell^{2}(\pi)$.

In what follows, the solution $u^{{\scriptscriptstyle 0}}=(u_n^{{\scriptscriptstyle 0}})_{n\in\entiers}$ of the homogeneous equation  
\begin{equation}\label{homog-eqn}
\lambda_{n}u_{n+1}+\mu_{n} u_{n-1}\un_{\{n\geq 2\}}-(\lambda_{n}+\mu_{n})u_{n}=0
\end{equation}
such that $u_1^{{\scriptscriptstyle 0}}=1$ will play an important role. Using \eqref{tricky} it is easy to verify that 
\begin{equation}\label{def:uno}
u^{{\scriptscriptstyle 0}}_{n}=1+\mathlarger{\sum}_{j=1}^{n-1}\frac{1}{\lambda_{j}\pi_{j}}, \; n\geq 1
\end{equation}
with the convention that $\sum_{j=1}^{{\scriptscriptstyle 0}}=0$. 

\begin{remark}
\label{casse}
Notice that $u^{{\scriptscriptstyle 0}}\notin \ell^2(\pi)$. Indeed, using \eqref{tricky}, observe that 
\[
u_n^{{\scriptscriptstyle 0}}\geq \frac{1}{\lambda_{n-1}\pi_{n-1}}=\frac{1}{\mu_n\pi_n}.
\]
Hence
\[
\sum_{n=1}^N \big(u_n^{{\scriptscriptstyle 0}}\big)^2 \pi_n \geq \sum_{n=1}^N \frac{1}{\mu_n^2\pi_n}.
\]
But by \eqref{eq:pimu} the last sum tends to $+\infty$ when $N$ goes to infinity.
\end{remark}

\subsection{Estimates of the largest eigenvalue and of the associated eigenvector}

Our first main result gives the behavior of $\grho_{0}$ and $\varphi$ as functions of $K$ when $K$ gets large. 
Recall that $x_*$ and $\ns$ are defined in \eqref{eq:xstar} and \eqref{eq:nstar}, respectively, and
that $u^{{\scriptscriptstyle 0}}=(u_n^{{\scriptscriptstyle 0}})_n$ is the solution of the homogeneous equation \eqref{homog-eqn}. 
The function $H$ is defined in \eqref{def:H} and recall that $H''(x_*)>0$ (see \eqref{eq:generic}).

\begin{theorem}[Estimates of $\grho_0$ and $\varphi$]
\label{thm:rho0-phi}
\leavevmode\\
For all $K>1$, we have
\begin{align*}
\grho_0(K) & =  \frac{\left(\sqrt{\frac{\lambda_1}{\mu_1}}-\sqrt{\frac{\mu_1}{\lambda_1}}\right)
\sqrt{K H''(x_*)}\,x_*\tilde \lambda(x_*)}{\sqrt{2\pi }}\ 
\mathlarger{\mathlarger{e}}^{-K\mathlarger{\int}_{0}^{x_*} \log \frac{\tilde{\lambda}(x)}{\tilde{\mu}(x)}\,\mathrm{d}x}\\
& \quad \times \mathsmaller{\left(1+ \mathcal{O}\left( \frac{(\log K)^3}{\sqrt{K}}\right)\right)}.
\end{align*}
Moreover, for all $K>1$,  we have
\[
\sup_{n\in\entiers}\left|\varphi_n(K) - V_n(K)\right| \leq \mathcal{O}(1)\grho_{0}(K)K \log K
\]
where
\[
V_n(K)=
\begin{cases}
u_n^{{\scriptscriptstyle 0}} & \textup{if}\; n\leq \ns\\
u_{\ns}^{{\scriptscriptstyle 0}} & \textup{if}\; n\geq \ns.
\end{cases}
\]
\end{theorem}

\tpg{proof-thm:rho0-phi}. 
Notice that the constant $c$ defined by
\begin{equation}
\label{lec}
c=-\int_{0}^{x_{*}}\log\frac{\tilde\mu(x)}{\tilde\lambda(x)}\;\mathrm{d}x
\end{equation}
is strictly positive by the assumptions on the functions $\tilde\lambda$, $\tilde\mu$.
It will appear several times later on. 

\begin{remark}
In the logistic case, one finds
\[
\grho_0(K)=\frac{(\tilde{\lambda}-\tilde{\mu})^2} {\sqrt{2\pi\tilde{\mu}}}\;\mathsmaller{\sqrt{K}}
\ \mathlarger{\mathlarger{e}}^{- K\left(\tilde{\lambda}-\tilde{\mu}+\tilde{\mu}\log\frac{{\tilde{\mu}}}{\tilde\lambda}\right)}\,
\mathsmaller{\left( 1+ \mathcal{O}\left(\frac{(\log K)^3}{\sqrt{K}}\right)\right).}
\]
\end{remark}

The following theorem provides a \emph{lower bound} for the spectral gap.  

\begin{theorem}[Spectral gap]
\label{rho1}
\leavevmode\\
There exists a constant $d>0$ such that for all $K>1$
\[
\grho_{1}(K) - \grho_{0}(K) \geq \frac{d}{\log K}.
\]
\end{theorem}

\tpg{proof-rho1}.

\begin{remark}
As a consequence of the preceding two theorems, one has
$\grho_{0}(K) \ll \grho_{1}(K)-\grho_{0}(K)$ for large $K$ because
$\grho_0(K)\approx \sqrt{K} \mathlarger{\mathlarger{e}}^{-cK}$.
\end{remark}

\subsection{Quasi-stationary distribution, survival rate and mean time to extinction}

We refer to \cite{MV} and \cite{CMSM} for background and more informations about quasi-stationary
distributions. As usual, we shall denote by $\proba_\nu$ the law of the process starting from a distribution
$\nu$ and by $\proba_n$ the law of the process starting from the state $n$, {\em i.e.} starting from the distribution $\delta_n$. The corresponding exepectations are respectively denoted by $\esperance_\nu$ and
$\esperance_n$.

\begin{proposition}\label{prop:qsd}
\leavevmode\\
For all $K>1$, the probability measure $\nu=(\nu_n)_n$ on $\entiers$ defined  by
\[
\nu_n=\frac{\pi_n \varphi_n}{\langle \varphi, \un\rangle_\pi}
\]
is the unique quasi-stationary distribution of the birth and death process.
\end{proposition}

Note that the quasi-stationary distribution $\nu$ depends on $K$ through $\varphi$.

\begin{proof}
In order to prove that $\nu$ is a quasi-stationary distribution, we must verify that
$\proba_{\nu}(X_t^{\scriptscriptstyle{K}} \in A | T_0>t)=\nu(A)$ for all $t>0$ and
for all subsets $A\subseteq\entiers$. 
Observe that for all $A\subseteq\entiers$, $\un_A\in\ell^2(\pi)$.
We have, using that $L$ is self-adjoint,
\begin{align*}
\proba_{\nu}(X_t^{\scriptscriptstyle{K}} \in A ,  T_0>t)
& =\mathlarger{\sum}_{n\in \entiers} \nu_{n} \,P_t(n,A)=
\frac{\langle \varphi, \mathlarger{\mathlarger{e}}^{tL}\un_A\rangle_\pi }{\langle \varphi, \un\rangle_\pi}\\
& =\frac{\langle \mathlarger{\mathlarger{e}}^{tL}\varphi, \un_A\rangle_\pi}{\langle \varphi, \un\rangle_\pi}
= \mathlarger{e}^{-\rho_0 t} \frac{\langle \varphi, \un_A\rangle_\pi}{\langle \varphi, \un\rangle_\pi}\\
& =  {\mathlarger{e}}^{-\rho_0 t}  \, \nu(A).
\end{align*}
Replacing $A$ by $\entiers$ yields the wanted relation. Since we have uniqueness (by \eqref{QSD}),
$\nu$ must be the quasi-stationary distribution.
\end{proof}

Before proceeding with the other results, we observe that the previous proof shows that for all $t>0$
\[
\proba_{\nu}(T_0>t)=\mathlarger{e}^{-\rho_0 t}.
\]
The quantity $\grho_0$ is usually called the exponential rate of survival. 
The mean time to extinction (starting from the quasi-stationary distribution) is thus 
\[
\esperance_{\nu}\big[T_0\big]=\frac{1}{\grho_0{\scriptstyle (K)}}.
\]

In view of Theorem \ref{thm:rho0-phi},  it is of order $e^{cK}/\sqrt{K}$ for some positive constant $c$. More precisely, we have the following corollary.

\begin{corollary}[Approximation of the mean time to extinction]
\leavevmode\\
For all $K>1$ we have
\begin{align*}
\esperance_{\nu}\big[T_0\big] & =  \frac{\sqrt{2\pi }}{\left(\sqrt{\frac{\lambda_1}{\mu_1}}-\sqrt{\frac{\mu_1}{\lambda_1}}\,\right)
\sqrt{K H''(x_*)}\,x_*\tilde \lambda(x_*)}\ 
\mathlarger{\mathlarger{e}}^{K\mathlarger{\int}_{0}^{x_*} \log \frac{\tilde{\lambda}(x)}{\tilde{\mu}(x)}\,\mathrm{d}x}\\
& \quad \times \mathsmaller{\left(1+ \mathcal{O}\left( \frac{(\log K)^3}{\sqrt{K}}\right)\right)}.
\end{align*}
\end{corollary}
Note that there is another way to obtain the above estimate of $\esperance_{\nu}\big[T_0\big]$. Indeed, we have
\[
\esperance_{\nu}\big[T_0\big]=\sum_{n\in\entiers}\esperance_{n}\big[T_0\big]\, \nu_n
\]
and since (see \cite{KT})
\[
\esperance_{n}\big[T_0\big]=\sum_{m=1}^n \frac{1}{\lambda_m \mu_m} \sum_{i\geq m+1} \pi_i,
\]
the estimate can be obtained by using Proposition \ref{prop:qsd} and Theorem \ref{thm:rho0-phi} to deal with $\varphi_n$.

\subsection{Convergence rate to the quasi-stationary distribution and Gaussian approximation}

We denote by $\mathlarger{\mathrm{d}}_{\scriptscriptstyle{\textup{TV}}}(\mu^{{\scriptscriptstyle (1)}},\mu^{{\scriptscriptstyle (2)}})$ the total variation distance between two probability measures $\mu^{{\scriptscriptstyle (1)}}$ and $\mu^{{\scriptscriptstyle (2)}}$. Recall that
\[
\mathlarger{\mathrm{d}}_{\scriptscriptstyle{\textup{TV}}}\big(\mu^{{\scriptscriptstyle (1)}},\mu^{{\scriptscriptstyle (2)}}\big)= 
\sup_{A\in \mathscr{P}(\integers)} \big|\mu^{{\scriptscriptstyle (1)}}(A)-\mu^{{\scriptscriptstyle (2)}}(A)\big|=\frac12\,
\mathlarger{\sum}_{n\in\integers}\, \big|\mu^{{\scriptscriptstyle (1)}}_n-\mu^{{\scriptscriptstyle (2)}}_n\big|
\]
where $\mathscr{P}(\integers)$ is the powerset of $\integers$. 

The process $\big(X^{\scriptscriptstyle{K}}_{t}, t\geq 0)$ is said to have a Yaglom limit if there exists a probability measure $\mathfrak m$ on $\entiers$ such that for every $n\in\entiers$ and for every $A\in \mathscr{P}(\entiers)$ one has
\[
\lim_{t\to\infty} \proba_n\left(X_t^{\scriptscriptstyle{K}} \in A \big| T_0>t\right)=\mathfrak m(A).
\]
When it exists, the Yaglom limit is a quasi-stationary distribution (whereas the converse is false in general), see \cite{MV}.

The following theorem provides a quantitative bound for the distance (in total variation) between 
the law of the process and a convex combination of the Dirac mass at $0$ and the quasi-stationary distribution
$\nu$. It also shows that $\nu$ is the Yaglom limit of $\big(X^{\scriptscriptstyle{K}}_{t}, t\geq 0)$ with a quantitative error bound. Recall that $-\grho_1$ is the second largest eigenvalue of $L$ (see Theorem \ref{autoadj}).

\begin{theorem}\label{thm:dtv}
\leavevmode\\
There exist three strictly positive constants $a, c_1,C$ such that for all $K>1$, for all $n\in\entiers$ and for all $t\geq 0$,
we have
\begin{align}
\nonumber
& \mathlarger{\mathrm{d}}_{\scriptscriptstyle{\textup{TV}}}
\left(\proba_{n}\big(X^{\scriptscriptstyle{K}}_{t}\in\,\cdot\,\big)\,,\,
\alpha_{n}{\scriptstyle(K)}\,\nu +(1-\alpha_{n}{\scriptstyle(K)})\gdelta_{0}\right)\\
& 
\le C\left({K}^{3/2}\log K \,\mathlarger{\mathlarger{e}}^{-c K}+
\big(1-\mathlarger{\mathlarger{e}}^{-\rho_{0}\,t}\big)+
K \mathlarger{\mathlarger{e}}^{-\frac{a}{4} t} + K^{\frac{3}{4}}\,
\mathlarger{\mathlarger{e}}^{c_1 K} \mathlarger{\mathlarger{e}}^{-\frac{\rho_1}{2} t}
\right)
\label{eq:convcomb}
\end{align}
where
\[
\alpha_{n}(K)=
\begin{cases}
\frac{u^{{\scriptscriptstyle 0}}_{n}}{u^{{\scriptscriptstyle 0}}_{n^*{\scriptscriptstyle (K)}}}\;&\mathrm{for}\;n\le \ns\\
1\;&\mathrm{for}\;n\ge \ns\\
\end{cases}
\]
and where $u^{{\scriptscriptstyle 0}}$ is defined in \eqref{def:uno}. Moreover 
\begin{equation}\label{eq:convyag}
\mathlarger{\mathrm{d}}_{\scriptscriptstyle{\textup{TV}}}
\left(\frac{P_t(n,\cdot)}{P_t(n,\entiers)}\, ,\, \nu\right)\leq 
C \left( K\, \mathlarger{\mathlarger{e}}^{-(\frac{a}{4}-\rho_0)t} + K^{\frac{3}{4}}\, \mathlarger{\mathlarger{e}}^{c_1 K}
\mathlarger{\mathlarger{e}}^{-(\frac{\rho_{1}}{2}-\rho_0)t}\right).
\end{equation}
In particular, the probability measure $\nu
$ is the Yaglom limit (in total variation distance) of the process $(X_t^{\scriptscriptstyle{K}}, t\geq 0)$.
\end{theorem}

\tpg{proof-thm:dtv}.

\begin{remark}
The proof of the previous theorem consists in establishing the following more explicit estimate:
there exist three strictly positive constants $a,c_1,C$ such that for all $K>1$, for all $n\in\entiers$ and for all $t\geq 0$,
we have
\begin{align}
\nonumber
& \mathlarger{\mathrm{d}}_{\scriptscriptstyle{\textup{TV}}}
\left(\proba_{n}\big( X^{\scriptscriptstyle{K}}_{t}\in\,\cdot\,\big)\,,\,
\frac{\langle \varphi,\un\rangle_{\pi}}{\|\varphi\|_{\pi}^{2}} \, \mathlarger{\mathlarger{e}}^{-\rho_0 t}\varphi_{n}\,
\nu(\cdot)+\Big(1-\frac{\langle \varphi,\un\rangle_{\pi}}{\|\varphi\|_{\pi}^{2}}\, \mathlarger{\mathlarger{e}}^{-\rho_0 t}\varphi_{n}\Big)\,\gdelta_{0}(\cdot)\right)\\
&\leq C\, \Big( K\, \mathlarger{\mathlarger{e}}^{-\frac{a}{4} t} + K^{\frac{3}{4}}\,
\mathlarger{\mathlarger{e}}^{c_1 K}\mathlarger{\mathlarger{e}}^{-\frac{\rho_1}{2} t}\Big).
\label{eq:maindtv}
\end{align} 
Then we show that the estimates \eqref{eq:convcomb} and \eqref{eq:convyag} follow from \eqref{eq:maindtv}.
\end{remark}

\begin{remark}
The estimate \eqref{eq:convcomb} can be interpreted as follows.
Recall that, for $K$ large, $\grho_{0}$ is very small. 
Therefore,  if we start with $n=\mathcal{O}(K)$ and if $t$ is such that
$K\log K/(\rho_{1}-\rho_{0})\ll t\ll 1/\grho_{0}$, we get the following rough estimate:
\[
\mathlarger{\mathrm{d}}_{\scriptscriptstyle{\textup{TV}}}\left(\proba_{n}\big(X^{\scriptscriptstyle{K}}_{t}\in\,\cdot\,\big)\,,\,
\alpha_{n}{\scriptstyle(K)}\nu +(1-\alpha_{n}{\scriptstyle(K)})\;\gdelta_{0}\right)\ll 1\;.
\]
This inequality highlights the existence of an interval of time during which the process is either extinct with a probability
close to $1-\alpha_{n}{\scriptstyle(K)}$ or obeys the quasi-stationary distribution $\nu$ with a probability close $\alpha_{n}{\scriptstyle(K)}$. 
This interval has a length that is roughly exponentially large in $K$.
\end{remark}

\begin{remark}
It follows from Theorem \ref{thm:rho0-phi} and Theorem \ref{rho1} that, for $K$ large enough,
\[
\min\left(\frac{a}{4}-\grho_0,\frac{\grho_{1}}{2}-\grho_0\right)\geq \frac{d}{3\log K}.
\]
Hence, for $K$ large enough,  the estimate \eqref{eq:convyag} can be written as
\[
\mathlarger{\mathrm{d}}_{\scriptscriptstyle{\textup{TV}}}
\left(\frac{P_t(n,\cdot)}{P_t(n,\entiers)}\, ,\, \nu\right)\leq 
2 C  K \mathlarger{\mathlarger{e}}^{c_1 K} \mathlarger{\mathlarger{e}}^{-\frac{d}{3\log K}t}.
\]
\end{remark}

\begin{remark}
Note that for every $n\geq 1$, the weights $\alpha_n{\scriptstyle(K)}$ appearing in \eqref{eq:convcomb}
can be written as
\[
\alpha_{n}{\scriptstyle(K)}=1-\left(\frac{\mu_{1}}{\lambda_{1}}\right)^{n}
+\frac{\Oun}{K}
\]
for all $K>1$. 
This follows by adapting the proof of Lemma \ref{b1}.
\end{remark}

The last result shows that the quasi-stationary distribution $\nu$ is close, as $K$ gets large,
to a Gaussian law centered at $\ns$. Recall that the function $H$ is defined in \eqref{def:H}.

\begin{theorem}\label{thm:ladqs}
\leavevmode\\
We have 
\[
\mathlarger{\mathrm{d}}_{\scriptscriptstyle{\textup{TV}}}\big(\nu^{\scriptscriptstyle{K}},G^{\scriptscriptstyle K}\big)
\le \frac{\Oun}{\sqrt K}
\]
where $G^{\scriptscriptstyle K}$ is the probability measure on $\entiers$ given by
\[
G^{\scriptscriptstyle K}_n=\frac{1}{Z(K)}\;
\mathlarger{\mathlarger{e}}^{-\frac{(n-n_{*}{\scriptscriptstyle (K)})^{2}}{2K\sigma^{2}}}
\]
where
\[
Z(K)=\sum_{n=1}^{\infty}
\mathlarger{\mathlarger{e}}^{-\frac{(n-n_{*}{\scriptscriptstyle (K)})^{2}}{2K\sigma^{2}}}
=\sqrt{2\pi K}\, \sigma+\Oun\;
\]
and where 
\[
\sigma=\frac{1}{\sqrt{H''(x_*)}}.
\]
\end{theorem}

Recall that $H''(x_*)>0$ by \eqref{Hseconde}. 
In the logistic case, one has $\sigma=\sqrt{\tilde{\lambda}}$.
\tpg{proof-thm:ladqs}.

\section{Proof of Theorem \ref{autoadj}}\label{proof-autoadj}

\subsection{$\geno$ is symmetric and closable in $\ell^{2}(\pi)$}

Using \eqref{tricky}, the reader can verify that,  for all $u,v\in \mathscr{D}$, one has
$\langle \tilde \gen u, v\rangle_\pi=\langle u, \tilde \gen v\rangle_\pi$. Hence $\tilde \gen$ is symmetric.\\
To verify closedness, one can apply a result in \cite[III.5.3]{Kato} saying that it is equivalent to
prove that, for every sequence $(y^{(k)})_{k}\in \mathscr{D}$ such that $y^{(k)}\to0$ (in $\ell^{2}(\pi)$) and 
such that $\tilde\gen y^{(k)}$ converges to $y$ (in $\ell^{2}(\pi)$), $y=0$. Details are left to the reader.

\subsection{$\gen$ defines a positive contraction semigroup in $\ell^{2}(\pi)$}

The key result in proving this claim is the following.
\begin{proposition}\label{resol}
\leavevmode\\
For every $f\in \ell^2(\pi)$ and every $\grho>0$, the equation 
\[
\big(\grho-\gen)y=f
\]
has a unique solution $y\in \mathcal{D}$ denoted by $\,R_{\rho} f$. Moreover
\[
\|R_{\rho} f\|_{\pi}\le \grho^{-1}\|f\|_{\pi}.
\]
Finally, if $f$ is nonnegative, so is $\,R_{\rho} f$. 
\end{proposition}

It is well-known that the previous bound is a sufficient condition for $\gen$ 
to generate a $C_{0}$ contraction semigroup $Q_{t}$ in $\ell^{2}(\pi)$, see {\em e.g.} \cite[p. 249]{yosida}.

The proof of this proposition requires two preliminary results. 
For $1\le n\le N$ we define (on $\ell^{{\scriptscriptstyle \infty}}{\scriptstyle(\{1,\ldots, N\})}$)
the truncated operator $L_{N}$ by
\[
(L_{{\scriptscriptstyle N}}v)_{n}=\lambda_{n}v_{n+1}\un_{\{n<N\}}+\mu_{n}v_{n-1}\un_{\{n\geq 2\}}
-\big(\lambda_{n}+\mu_{n}\big)v_n.
\]
The operator $L_{{\scriptscriptstyle N}}$ satisfies the following positive maximum principle.
\begin{lemma}\label{pm}
\leavevmode\\
Let $v\in \ell^{{\scriptscriptstyle \infty}}{\scriptstyle(\{1,\ldots,N\})}$ and let $m\in \{1,\ldots,N\}$
such that $v_m=\sup_{1\le n\le N}v_n$. \\ If $v_m\ge 0$, then $(L_{{\scriptscriptstyle N}}v)_m\le 0$.
\end{lemma}

\begin{proof}
For $2\leq m\leq N-1$, we get
\[
(L_{{\scriptscriptstyle N}}v)_m =
\lambda_{m} v_{m+1} + \mu_{m} v_{m-1} -( \lambda_{m} +\mu_{m})v_{m}
\leq 0
\]
since, by definition of  $m$, $v_m$ is maximal. 
The cases $m=1$ and $m=N$ follow similarly.
\end{proof}

\begin{lemma}\label{conseqpm}
\leavevmode\\
Let  $g\in \ell^{{\scriptscriptstyle \infty}}{\scriptstyle(\{1,\ldots,N\})}$ and $\grho>0$.
The equation $(\grho-L_{{\scriptscriptstyle N}})v=g$ has a unique solution in  $\ell^{{\scriptscriptstyle \infty}}{\scriptstyle(\{1,\ldots,N\})}$.
Moreover,  one has $\|v\|_{\ell^{\infty}{\scriptscriptstyle(\{1,\ldots,N\})}}
\le \|g\|_{\ell^{\infty}{\scriptscriptstyle(\{1,\ldots,N\})}}/\grho$. Finally, if $g\ge0$ then $v\ge 0$.
\end{lemma}

\begin{proof}
If $g\in \ell^{{\scriptscriptstyle \infty}}{\scriptstyle(\{1,\ldots,N\})}$ and $\grho>0$ are such that
$g=(\grho-L_{{\scriptscriptstyle N}})v$, and if $m\in\{1,\ldots,N\}$ is such that
$v_{m}=\sup_{1\le n\le N}v_n\ge 0$
then, by Lemma \ref{pm}, $v_{m}\le g_{m}/\grho$.
Considering $-v$ and $-g$, it follows that if $(\grho-L_{{\scriptscriptstyle N}})v=g$ and if
$m\in\{1,\ldots,N\}$ is such that
$v_{m}=\inf_{1\le n\le N}v_n\le 0$ then $v_{m}\ge g_{m}/\grho$.
This implies that $v\geq 0$ if $g\geq 0$.
The previous two inequalities imply
$\|v\|_{\ell^{\infty}{\scriptscriptstyle(\{1,\ldots,N\})}} \le \|g\|_{\ell^{\infty}{\scriptscriptstyle(\{1,\ldots,N\})}}/\grho$.
In particular we have $\textup{Ker}(\grho-L_{{\scriptscriptstyle N}})=\{0\}$, namely $\grho-L_{{\scriptscriptstyle N}}$ is
invertible in $\ell^{{\scriptscriptstyle \infty}}{\scriptstyle(\{1,\ldots,N\})}$.
The lemma is proved.
\end{proof}

We now turn to the proof of Proposition \ref{resol}.\\
Let $f\in \mathscr{D}$ and let $N_0\geq 1$ be such that $f_n=0$ for all $n>N_0$.
Applying Lemma \ref{conseqpm} for $N>N_0$ yields a $v^{\scriptscriptstyle{(N)}}\in\ell^{{\scriptscriptstyle \infty}}{\scriptstyle(\{1,\ldots,N\})}$
such that
\[
(\grho-L_{{\scriptscriptstyle N}})\;v^{\scriptscriptstyle{(N)}}=f\;.
\]
We also have that for all $N>N_0$
\begin{equation}\label{eq:vN}
\|v^{\scriptscriptstyle{(N)}}\|_{\ell^{\infty}{\scriptscriptstyle(\{1,\ldots,N\})}}\leq \frac{1}{\grho}\|f\|_{\ell^{\infty}{\scriptscriptstyle(\{1,\ldots,N_0\})}}.
\end{equation}
Define $u^{\scriptscriptstyle{(N)}}\in\mathscr{D}$ by
\[
u^{\scriptscriptstyle{(N)}}_n=
\begin{cases}
v^{\scriptscriptstyle{(N)}}_n & \mathrm{if}\quad n\le N\\
0 & \mathrm{if}\quad n>N.
\end{cases}
\]
For all $p\in\entiers$ we have
\[
\big((\grho-\gen)\;u^{\scriptscriptstyle{(N)}}\big)_p=
f_p+\big[(\grho+\lambda_{{\scriptscriptstyle N}}+\mu_{{\scriptscriptstyle N}})v^{\scriptscriptstyle{(N)}}_{{\scriptscriptstyle N}}-
\mu_{{\scriptscriptstyle N}} v^{\scriptscriptstyle{(N)}}_{{\scriptscriptstyle N-1}}\big]\un_{\{p={\scriptscriptstyle N}\}}
-\mu_{{\scriptscriptstyle N+1}}v^{\scriptscriptstyle{(N)}}_{{\scriptscriptstyle N}}\un_{\{p={\scriptscriptstyle N+1}\}}.
\]
It is then easy to show that 
\[
\big\|(\grho-\gen)\;u^{\scriptscriptstyle{(N)}} - f\big\|_{\pi}^2 \leq \mathcal{O}(1) 
(\pi_{{\scriptscriptstyle N}}\mu_{{\scriptscriptstyle N}}^2+\pi_{{\scriptscriptstyle N+1}}\mu_{{\scriptscriptstyle N+1}}^2)
\]
by using \eqref{infini} and \eqref{eq:vN}.
Hence, since we assume that \eqref{eq:pimu} holds, we get that $(\grho-\gen)\;u^{\scriptscriptstyle{(N)}}$
converges strongly to $f$.
Using $u^{\scriptscriptstyle{(N)}}_{{\scriptscriptstyle N+1}}=0$ we obtain
\begin{equation}\label{eq:urhpLU}
\langle u^{\scriptscriptstyle{(N)}}, (\grho-\gen)u^{\scriptscriptstyle{(N)}}\rangle_{\pi}=
\langle u^{\scriptscriptstyle{(N)}},f\rangle_{\pi}+r_{\!{\scriptscriptstyle N}}
\end{equation}
where $r_{\!{\scriptscriptstyle N}}=\big[(\grho+\lambda_{{\scriptscriptstyle N}}+\mu_{{\scriptscriptstyle N}})
v^{\scriptscriptstyle{(N)}}_{{\scriptscriptstyle N}}-\mu_{{\scriptscriptstyle N}}
v^{\scriptscriptstyle{(N)}}_{{\scriptscriptstyle N-1}}\big]v^{\scriptscriptstyle{(N)}}_{{\scriptscriptstyle N}} \pi_{{\scriptscriptstyle N}}$.

One gets (recall that $u^{\scriptscriptstyle{(N)}}_{{\scriptscriptstyle N+1}}=0$)
\begin{align*}
& \langle u^{\scriptscriptstyle{(N)}},\gen u^{\scriptscriptstyle{(N)}}\rangle_{\pi}  \\
& =\sum_{n=1}^{\infty} \pi_n \lambda_{n} u^{\scriptscriptstyle{(N)}}_n u^{\scriptscriptstyle{(N)}}_{n+1}+
\sum_{n=2}^{\infty} \pi_n\mu_{n} u^{\scriptscriptstyle{(N)}}_n u^{\scriptscriptstyle{(N)}}_{n-1}
-\sum_{n=1}^{\infty} \pi_n\big(\lambda_{n}+\mu_{n}\big) (u^{\scriptscriptstyle{(N)}}_{n})^{2}
\end{align*}
\begin{align*}
& \leq \frac12 \sum_{n=1}^{\infty} \pi_n \lambda_{n}(u^{\scriptscriptstyle{(N)}}_{n})^2
+
\frac12 \sum_{n=1}^{\infty} \pi_n \lambda_{n}(u^{\scriptscriptstyle{(N)}}_{n+1})^2
+
\frac12 \sum_{n=2}^{\infty} \pi_n \mu_{n}(u^{\scriptscriptstyle{(N)}}_{n})^2
\\
& \quad \;+
\frac12 \sum_{n=2}^{\infty} \pi_n \mu_{n} (u^{\scriptscriptstyle{(N)}}_{n-1})^2-\sum_{n=1}^{\infty} \pi_n\big(\lambda_{n}+\mu_{n}\big)(u^{\scriptscriptstyle{(N)}}_{n})^{2}\\
& \leq \frac12 \sum_{n=1}^{\infty} \pi_n \lambda_{n}\;(u^{\scriptscriptstyle{(N)}}_{n})^2
+
\frac12 \sum_{n=1}^{\infty} \pi_{n+1}\mu_{n+1}(u^{\scriptscriptstyle{(N)}}_{n+1})^2
+
\frac12 \sum_{n=2}^{\infty} \pi_n \mu_{n}(u^{\scriptscriptstyle{(N)}}_{n})^2
\\
& \quad \; +
\frac12 \sum_{n=1}^{\infty} \pi_n \lambda_{n}(u^{\scriptscriptstyle{(N)}}_{n})^2-
\sum_{n=1}^{\infty} \pi_n\big(\lambda_{n}+\mu_{n}\big)(u^{\scriptscriptstyle{(N)}}_{n})^{2}
\\
& \leq -\pi_1 \mu_1  (u^{\scriptscriptstyle{(N)}}_{1})^2\leq 0
\end{align*}
where we used \eqref{tricky}.
Hence, it follows from \eqref{eq:urhpLU} and the previous inequality that 
\begin{align*}
\grho \|u^{\scriptscriptstyle{(N)}}\|^{2}_{\pi} & =\langle u^{\scriptscriptstyle{(N)}},(\rho-L)u^{\scriptscriptstyle{(N)}}\rangle_{\pi}
+ \langle u^{\scriptscriptstyle{(N)}},Lu^{\scriptscriptstyle{(N)}}\rangle_{\pi}\\
& \le \langle u^{\scriptscriptstyle{(N)}},f\rangle_{\pi}+r_{\!{\scriptscriptstyle N}}
\le \|u^{\scriptscriptstyle{(N)}}\|_{\pi}\,\|f\|_{\pi}+r_{\!{\scriptscriptstyle N}}.
\end{align*}
Therefore we obtain
\[
\|u^{\scriptscriptstyle{(N)}}\|_{\pi}\le \frac{\|f\|_{\pi}}{2\grho} +\sqrt{\frac{r_{{\scriptscriptstyle N}}}{\grho} + \frac{\|f\|_{\pi}^2}{4\grho^2}}
\]
where the right hand side is the largest root of the polynomial function
$x\mapsto \rho x^2-\|f\|_\pi x-r\!_{{\scriptscriptstyle N}}$.
Since $r_{{\scriptscriptstyle N}}$ tends to $0$ by \eqref{eq:pimu} when $N$ tends to infinity,
$\sup_{{\scriptscriptstyle N}}\|u^{\scriptscriptstyle{(N)}}\|_{\pi}<\infty$.
Since a ball in the Hilbert space $\ell^{2}(\pi)$ is weakly compact \cite[p. 126]{yosida}, we can extract from
the sequence $(u^{\scriptscriptstyle{(N)}})$ a subsequence weakly converging to some $u\in\ell^{2}(\pi)$.
Moreover
\[
\|u\|_{\pi}\le \frac{1}{\grho}\|f\|_{\pi}
\]
by \cite[Theorem 1, p. 120]{yosida}.
Since the sequence $((\grho-\gen)u^{{\scriptscriptstyle (N)}})$ is also weakly convergent to $f$ (see above, even strongly convergent in our case),
we can apply \cite[Problem 5.12, p. 165]{Kato} to conclude that $u\in\mathcal{D}$ and $(\grho-\gen)u=f$.\\
At this point, we have proved that for all $f\in\mathscr{D}$ the equation $(\grho-\gen)u=f$ has a solution in $\mathcal{D}$.

If $f$ is nonnegative, Lemma \ref{conseqpm} implies that all the
$u^{\scriptscriptstyle{(N)}}$ are nonnegative for $N$ large enough, hence $u$ is
nonnegative. 

For every $w\in\mathcal{D}$, there is a sequence $(w^{(n)})$, with $w^{(n)}\in\mathscr{D}$ for all $n$, converging to 
$w$ (in $\ell^{2}(\pi)$) with $(\gen w^{(n)})$ converging to $\gen w$ in
$\ell^{2}(\pi)$ (see \cite[III.5.2]{Kato}). As before, 
\[
\langle w^{(n)},(\grho-\gen)w^{(n)}\rangle_{\pi} \ge \grho\|w^{(n)}\|^{2}_{\pi}.
\]
Therefore
\begin{equation}\label{Lneg}
\langle w,(\grho-\gen)w\rangle_{\pi} \ge \grho\|w\|^{2}_{\pi}
\end{equation}
for all $w\in\mathcal{D}$.
This implies that the equation
\[
(\grho-\gen)u=f.
\]
has a unique solution $u\in\mathcal{D}$ for every $f\in\mathscr{D}$.
This solution, denoted by $R_{\rho}f$, satisfies
\[
\big\|R_{\rho}f\big\|_{\pi}\le \frac{1}{\grho} \|f\|_{\pi}
\]
and it is nonnegative if $f$ is nonnegative. Since $\mathscr{D}$ is dense
in $\ell^{2}(\pi)$, the linear operator $R_{\rho}$ can be extended to a linear
operator on $\ell^{2}(\pi)$ with a norm that is at most $\grho^{-1}$ (see \cite[II.2.2]{Kato}).

Since $\mathscr{D}$ is dense in $\ell^{2}(\pi)$, for each $f\in \ell^{2}(\pi)$ we can find a sequence $(f^{(k)})\subset \mathscr{D}$ converging to $f$ in $\ell^{2}(\pi)$. Moreover, $\big(R_{\rho}f^{(k)})$ converges to $R_{\rho}f$.
Since, for all $k$,  $R_{\rho}f^{(k)}\in \mathcal{D}$ and 
$
\gen R_{\rho} f^{(k)}=\grho R_{\rho}f^{(k)}-f^{(k)}
$
converges in $\ell^2(\pi)$ to $\grho R_{\rho}f-f$, we conclude, by using \cite[III.5.2]{Kato},
that, for every $f\in\ell^2(\pi)$, $R_{\rho} f\in \mathcal{D}$ and
\[
(\grho-\gen) R_{\rho}f=f.
\]
Nonnegativity follows easily. This finishes the proof of the proposition.

We can now make the proof of statement 2 in Theorem \ref{autoadj}.
Using Proposition \ref{resol}, we can apply \cite[p. 249]{yosida} to show that $\gen$ generates a
$C_{0}$ contraction semigroup $Q_{t}$ in $\ell^{2}(\pi)$. 
For all $t\ge0$, the operator  $Q_{t}$  maps nonnegative sequences
to nonnegative sequences since this holds for $R_{\rho}$ for all
$\grho>0$ using \cite[formula 3, p. 246]{yosida}. 

\subsection{Compactness, self-adjointness and dissipativity}

\paragraph{$\gen$ has a compact resolvent in $\ell^{2}(\pi)$.}

From the equation $\ (\grho - L)\, R_{\rho}= \mathrm{Id}$, we get for every $f\in\ell^2(\pi)$
\[
(R_{\rho}f)_{n}=-\frac{f_{n}}{\grho+\lambda_{n}+\mu_{n}}+
\frac{\lambda_{n}(R_{\rho}f)_{n+1}}{\lambda_{n}+\mu_{n}+\grho}
+\frac{\mu_{n}(R_{\rho}f)_{n-1}}{\lambda_{n}+\mu_{n}+\grho}\ \un_{\{n\geq 2\}}.
\]
We are going to verify that each term is uniformly square summable at infinity with respect to
the weights $(\pi_{n})$. \\
This is obvious for the first term since $\lim_{n\to\infty}\frac{1}{\lambda_{n}+\mu_{n}+\grho}=0$.\\
For the other two terms, by using \eqref{tricky}, we have for all $N\geq 2$
\[
\sum_{n=N}^\infty  \left|(R_{\rho}f)_{n+1}\right|^2 \frac{\lambda_{n}^2\, \pi_{n}}{(\lambda_{n}+\mu_{n}+ \grho)^2} 
+
\sum_{n=N}^\infty  \left|(R_{\rho}f)_{n-1}\right|^2 \frac{\mu_{n}^2\, \pi_{n}}{(\lambda_{n}+\mu_{n}+ \grho)^2}
=
\]
\[
\sum_{n=N}^\infty  \left|(R_{\rho}f)_{n+1}\right|^2 \frac{\lambda_{n}\,\mu_{n+1} \pi_{n+1}}{(\lambda_{n}+\mu_{n}+ \grho)^2} 
+
\sum_{n=N}^\infty 
\left|(R_{\rho}f)_{n-1}\right|^2 \frac{\mu_{n}\,\lambda_{n-1}\, \pi_{n-1}}{(\lambda_{n}+\mu_{n}+ \grho)^2}.
\]
Using \eqref{infini} and \eqref{eq:mumu} we conclude that for all $\varepsilon>0$, there exists $N_{\varepsilon}$ such that for all $N\geq N_{\varepsilon}$
\begin{align*}
& \sum_{n=N}^\infty\!\! \left( \!\left|(R_{\rho}f)_{n+1}\right|^2 \frac{\lambda_{n}^2\, \pi_{n}}{(\lambda_{n}+\mu_{n}+ \grho)^2} 
+\!
\left|(R_{\rho}f)_{n-1}\right|^2 \frac{\mu_{n}^2\, \pi_{n}}{(\lambda_{n}+\mu_{n}+ \grho)^2}\!\right)\\
& \leq \varepsilon \|R_{\rho}f\|^2_{\pi} \leq \frac{\varepsilon \|f\|^2_{\pi}}{\rho^2}.
\end{align*}
Compactness of the resolvent follows.

If $-\grho$ is an eigenvalue, a corresponding eigenvector $u$ (in $\ell^2(\pi)$) must satisfy the identities 
\begin{align*}
& u_{2} = \frac{(\lambda_{1} + \mu_{1}-\grho)u_{1}}{\lambda_{1}},\\
& u_{n+1} = \frac{(\lambda_{n} + \mu_{n}-\grho)u_{n}}{\lambda_{n}}-
\frac{\mu_{n}u_{n-1}}{\lambda_{n}},\quad \forall n\geq 2.
\end{align*}
Therefore, $u_{1}$ determines all the ${u_{n}}'s$. This implies that all eigenvalues are simple.\\
Positivity of the eigenvector associated with the maximal eigenvalue $\,-\grho_{0}$ follows from
the fact that the semigroup preserves nonnegativity and the fact that if an eigenvector
is orthogonal to any positive function, it would be equal to $0$, which is not true.

\paragraph{Self-adjointness and dissipativity.}
Self-adjointness follows by an argument found in \cite[problem V.3.32, p. 279]{Kato}.
In more details, it follows from equation \eqref{Lneg} that for all $u\in \mathcal{D}$,
$\langle u,Lu\rangle_{\pi}\le 0$, hence is $L$ is dissipative
and the numerical range of $L$ is contained in the negative real
line. By Theorem V.3.2 page 268 in \cite{Kato} the defect index is
constant outside the negative real line, and equal to zero on the
positive real line by Proposition \ref{resol}.
Therefore the spectrum of $L$ is contained in the negative real line and $L$ is self adjoint
by Theorem 3.16 in \cite[Chapter V, p. 271]{Kato}.

\section{Proof of Theorem \ref{thm:rho0-phi}}\label{proof-thm:rho0-phi}

For every small number $\grho$, we are going to consider sequences $(u_{n})_{n}$ satisfying
\begin{equation}
\label{eqpropre}
\lambda_{n}u_{n+1}+\mu_{n} u_{n-1}\un_{\{n\geq 2\}}-(\lambda_{n}+\mu_{n})u_{n}= - \grho\, u_{n}.
\end{equation}
The strategy will be as follows.
If $\grho=0$, $u^{{\scriptscriptstyle 0}}$ is a solution of \eqref{eqpropre} for all $n\geq 1$ and the constant sequence $\,1\,$  is a solution of \eqref{eqpropre} for all $n\geq 2$. For small $\grho\neq 0$ and $n\le \ns$, we will look for a solution of \eqref{eqpropre} that is a small perturbation of $u^{{\scriptscriptstyle 0}}$.  Since $u^{{\scriptscriptstyle 0}} \notin \ell^2(\pi)$ (see Remark \ref{casse}), we cannot use such an argument for large $n$. For $n\ge \ns-1$, we will use Levinson's technique (see \cite{levinson}, \cite{fedor}) to prove that there is a solution of \eqref{eqpropre} that is almost constant.  Then we will match these two solutions in $\{\ns-1, \ns\}$.  This will be possible for a single value of $\grho$ that has to be $\grho_{0}$. Since \eqref{eqpropre}  is a recursion of order $2$, this matched sequence is a solution for all $n\in\entiers$. Finally we will prove that this sequence belongs to $\mathcal D$ (see Theorem \ref{autoadj} for the definition of $\mathcal D$).

\subsection{When $1\leq n \leq \ns$}

We look for a solution of the form
\[
v_{n}=u^{{\scriptscriptstyle 0}}_{n}\,(1+\delta_{n})
\]
where $u^{{\scriptscriptstyle 0}}=(u^{{\scriptscriptstyle 0}}_{n})$ is defined in \eqref{def:uno}.

\begin{proposition}\label{ledelta}
\leavevmode\\
There exists a constant $\widetilde{C}>0$ such that for $K$ large enough and for each
$\grho\in \left[-{\scriptstyle 1/(3\widetilde{C}K\log K)},{\scriptstyle 1/(3\widetilde{C}K\log K)}\right]$
the equation \eqref{eqpropre} admits for all $n\leq \ns$ a solution of the
form
\[
v_n=u^{{\scriptscriptstyle 0}}_{n}\,(1+\delta_{n})
\]
where
\begin{enumerate}
\item $\delta_1=0$;
\item $\delta_n$ is a solution of 
\[
\lambda_{n}\frac{u^{{\scriptscriptstyle 0}}_{n+1}}{u^{{\scriptscriptstyle 0}}_{n}}\big(\delta_{n+1}-\delta_{n}\big)-\mu_{n}
\frac{u^{{\scriptscriptstyle 0}}_{n-1}}{u^{{\scriptscriptstyle 0}}_{n}}
\big(\delta_{n}-\delta_{n-1}\big)\un_{\{n\geq 2\}}
=-\grho\big(1+\delta_{n}\big);
\]
\item  $1+\delta_{n} >0$ and 
$
\| (\delta_n)\|_{\ell^\infty{\scriptscriptstyle(\{1,\ldots, n_{*} (K)\})}}
\leq \frac{|\grho|\widetilde{C}K\log K}{1-|\grho|\widetilde{C}K\log K}. 
$
\item $\delta=(\delta_n)_n$ is a smooth function of $\grho$ and
\[
\left\|\frac{\mathrm{d}\delta}{\mathrm{d}{\grho}}(\grho)  - \Delta^{\!{\scriptscriptstyle 0}}\right\|_{\ell^\infty{\scriptscriptstyle(\{1,\ldots, n_{*} (K)\})}}
\leq 4(\widetilde{C}K\log K)^2\,|\grho|
\]
where
\begin{equation}\label{def:Deltan0}
\Delta^{\!{\scriptscriptstyle 0}}_n=
\sum_{j=1}^{n-1} \frac{1}{\lambda_j \pi_j u_j^{{\scriptscriptstyle 0}} u_{j+1}^{{\scriptscriptstyle 0}}}
+ \sum_{j=1}^{n-1} \sum_{p=2}^{j}\frac{(u_p^{{\scriptscriptstyle 0}})^2 \pi_p}{\lambda_j\pi_j u_j^{{\scriptscriptstyle 0}} u_{j+1}^{{\scriptscriptstyle 0}}} \quad \textup{for all}\;n\geq 2
\end{equation}
and $\Delta^{\!{\scriptscriptstyle 0}}_1=0$.
\end{enumerate}
\end{proposition}

\begin{proof}
It is easy to check that 
\[
\lambda_{n}\frac{u^{{\scriptscriptstyle 0}}_{n+1}}{u^{{\scriptscriptstyle 0}}_{n}}\,\big(\delta_{n+1}-\delta_{n}\big)
-\mu_{n}\frac{u^{{\scriptscriptstyle 0}}_{n-1}}{u^{{\scriptscriptstyle 0}}_{n}}\big(\delta_{n}-\delta_{n-1}\big)\un_{\{n\geq 2\}}
=-\grho\big(1+\delta_{n}\big).
\]
We impose $\delta_{1}=0$ ({\em i.e.} $v_1=1$).\\
We now apply Lemma \ref{soleqlin} for $n\ge2$ with
\[
h_{n}=-\grho\,(1+\delta_{n})\;,\quad  
\alpha_{n}=\lambda_{n}\,\frac{u^{{\scriptscriptstyle 0}}_{n+1}}{u^{{\scriptscriptstyle 0}}_{n}}\;,
\quad \beta_{n}=\mu_{n}\,
\frac{u^{{\scriptscriptstyle 0}}_{n-1}}{u^{{\scriptscriptstyle 0}}_{n}}.
\]
For $r>s$, we have
\[
\Theta_{r,s}=
\frac{\beta_{r-1}\ldots\beta_{s}}{\alpha_{r-1}\ldots\alpha_{s}}
=\frac{\mu_{r-1}\ldots\mu_{s}}{\lambda_{r-1}\ldots\lambda_{s}}\;
\frac{u^{{\scriptscriptstyle 0}}_{s-1}}{u^{{\scriptscriptstyle 0}}_{r-1}}\;\frac{u^{{\scriptscriptstyle 0}}_{s}}{u^{{\scriptscriptstyle 0}}_{r}}=
\frac{\lambda_{s-1}\pi_{s-1}}{\lambda_{r-1}\pi_{r-1}}\;
\frac{u^{{\scriptscriptstyle 0}}_{s-1}}{u^{{\scriptscriptstyle 0}}_{r-1}}\;\frac{u^{{\scriptscriptstyle 0}}_{s}}{u^{{\scriptscriptstyle 0}}_{r}}.
\]
Observing that $\lambda_1\, u_2^{{\scriptscriptstyle 0}}\,\delta_2=-\grho$, we get
\begin{equation}\label{deltan}
\delta_n= -\grho \left(\sum_{j=1}^{n-1} \frac{1}{\lambda_j \pi_j u_j^{{\scriptscriptstyle 0}} u_{j+1}^{{\scriptscriptstyle 0}}}
+ \sum_{j=1}^{n-1} \sum_{p=2}^{j}
 \frac{(u_p^{{\scriptscriptstyle 0}})^2 \pi_p(1+\delta_p)}{\lambda_j\pi_j u_j^{{\scriptscriptstyle 0}} u_{j+1}^{{\scriptscriptstyle 0}}}\right).
\end{equation}
Equation \eqref{deltan} can be written as
\[
\delta= -\grho \Delta^{\!{\scriptscriptstyle 0}} + \grho B \delta
\]
where $B$ is a linear operator defined as
\[
\big( B\delta\big)_n =
-\sum_{j=1}^{n-1} \sum_{p=2}^{j}\frac{(u_p^{{\scriptscriptstyle 0}})^2 \pi_p}{\lambda_j\pi_j u_j^{{\scriptscriptstyle 0}}
u_{j+1}^{{\scriptscriptstyle 0}}}\ \delta_p.
\]
Using Lemma \ref{propuno} and the fact that $\mu_\ell/\lambda_\ell<1$ for $\ell\leq \ns-1$, we have the bound
\begin{align*}
\Delta^{\!{\scriptscriptstyle 0}}_{\ns} & = \sum_{j=1}^{n_{*}{\scriptscriptstyle (K)}-1} \frac{1}{\lambda_j \pi_j
u_j^{{\scriptscriptstyle 0}} u_{j+1}^{{\scriptscriptstyle 0}}}
+ \sum_{j=1}^{n_{*}{\scriptscriptstyle (K)}-1} \sum_{p=2}^{j}
 \frac{(u_p^{{\scriptscriptstyle 0}})^2 \pi_p}{\lambda_j\pi_j u_j^{{\scriptscriptstyle 0}} u_{j+1}^{{\scriptscriptstyle 0}}}\\
&
\leq \Oun \left( \sum_{j=1}^{n_{*}{\scriptscriptstyle (K)}-1} \frac{1}{\lambda_j \pi_j}  +
\sum_{j=1}^{n_{*}{\scriptscriptstyle (K)}-1} \sum_{p=2}^{j}\frac{1}{\lambda_p}\right) \\
&
\leq \widetilde{C} K \log K
\end{align*}
where $\widetilde{C}>0$ is a constant independent of $K$ since $\lambda_p\geq p\tilde{\lambda}(0)$.
Therefore
\begin{align}
\notag
& \|\Delta^{\!{\scriptscriptstyle 0}}\|_{\ell^\infty{\scriptscriptstyle(\{1,\ldots, n_{*} (K)\})}}
\leq \widetilde{C} K\log K\\ 
& \textup{and}\quad
\|B\|_{\ell^\infty{\scriptscriptstyle(\{1,\ldots, n_{*} (K)\})}}
\leq  \Delta^{\!{\scriptscriptstyle 0}}_{n_{*}{\scriptscriptstyle (K)}} \le  \widetilde{C} K\log K.
\label{bornesbis}
\end{align}
We denote by $\Omega$ the complex disk centered at the origin and of radius
$
\frac{1}{3\widetilde{C} K\log K}.
$
For every $\grho\in\Omega$, the operator $\mathrm{Id} - \grho B$ is invertible and  $\delta = (\mathrm{Id} - \grho B)^{-1} \,\grho \Delta^{\!{\scriptscriptstyle 0}}$.
It follows from \eqref{bornesbis} that
\[
\|\delta\|_{\ell^\infty{\scriptscriptstyle(\{1,\ldots, n_{*} {\scriptscriptstyle (K)}\})}}\leq \frac{|\grho|\,\widetilde{C} K\log K}{1 -|\grho|\,\widetilde{C} K\log K}.
\]
Therefore, $\delta$ is bounded in $\ell^\infty{\scriptstyle(\{1,\ldots, n_{*} {\scriptscriptstyle (K)}\})}$ by $\frac{1}{2}$ and $1+\delta_{n}>0$ for all $n\leq \ns$. 
It also follows that $\delta=(\delta_n)_{1\leq n_{*}{\scriptscriptstyle (K)}}$ is an analytic function on $\Omega$.
We now compute its derivative in $\Omega$:
\[
\frac{\mathrm{d}\delta}{\mathrm{d}\grho} (\grho)=
(\mathrm{Id} - \grho B)^{-1} \, \Delta^{\!{\scriptscriptstyle 0}} + (\mathrm{Id} - \grho B)^{-2} \,\grho B\,\Delta^{\!{\scriptscriptstyle 0}}.
\]
Using \eqref{bornesbis} we get that for every $\grho\in\Omega$
\[
\left\|\frac{\mathrm{d}\delta}{\mathrm{d}\grho}(\grho)- \Delta^{\!{\scriptscriptstyle 0}}\right\|_{\ell^\infty{\scriptscriptstyle(\{1,\ldots, n_{*} (K)\})}} \le 4(\widetilde{C}K\log K)^2\, |\grho|.
\]
This finishes the proof of the proposition.
\end{proof}

\subsection{When $n\ge \ns-1$}

\begin{proposition}\label{lew}
\leavevmode\\
Let $C$ be the constant defined in Lemma \ref{controle}.
For $K$ large enough and each $\grho\in\left[-{\scriptstyle 1/(3CK)},{\scriptstyle 1/(3CK)}\right]$
the equation \eqref{eqpropre} admits for all $n\geq \ns-1$ a solution
\[
v_n = 1+w_n,
\]
where
\begin{enumerate}
\item $w_{\ns -1}=0$;
\item $w_n$ is a solution of $\lambda_{n}(w_{n+1}-w_{n})+\mu_{n}(w_{n-1}-w_{n})=- \grho(1+w_{n})$;
\item  $1+w_{n} >0$  and 
$
\| w_n\|_{\ell^\infty(\{n_*{\scriptscriptstyle (K)} -1,n_*{\scriptscriptstyle (K)} ,\ldots\})}\leq \frac{|\grho|CK}{1-|\grho|CK}.
$
\item $w=(w_n)$ is a smooth function of $\grho$ and
\[
\left\|\frac{\mathrm{d}w}{\mathrm{d}\grho}(\grho) -
W^{{\scriptscriptstyle 0}}\right\|_{\ell^\infty(\{n_*{\scriptscriptstyle (K)} -1,n_*{\scriptscriptstyle (K)} ,\ldots\})}
\le 4(CK)^2|\grho|
\]
where 
\begin{equation}\label{def-Wn0}
W_n^{{\scriptscriptstyle 0}} = \sum_{j=n_{*}{\scriptscriptstyle (K)}-1} ^{n-1}
\sum_{p=j+1}^{\infty}  \frac{ \pi_{p}}{\lambda_{j}\pi_{j}}\quad\textup{for all}\;n\geq \ns
\end{equation}
and $W_{n_*{\scriptscriptstyle (K)} -1}=0$.
\end{enumerate}
\end{proposition}

\begin{proof}
Let us define by induction for $n\geq \ns$,
\begin{equation}
\label{wn}
w_{n}=  \grho \sum_{j=n_*{\scriptscriptstyle (K)} -1} ^{n-1}\sum_{p=j+1}^{\infty}  \frac{ \pi_{p}}{\lambda_{j}\pi_{j}}(1+w_{j}),
\end{equation}
with $w_{n_{*}{\scriptscriptstyle (K)} -1}=0$.
It is easy to check by using \eqref{tricky} that 
\[
\lambda_{n}\,w_{n+1}+\mu_{n}\,w_{n-1}
-(\lambda_{n}+\mu_{n})w_{n}=-\grho (1+w_{n}).
\]
Equation \eqref{wn} can be written as
\[
w = \grho W^{{\scriptscriptstyle 0}} + \grho A w
\]
where $A$ is a linear operator defined as
\[
\big( A w\big)_n=\sum_{j=n_*{\scriptscriptstyle (K)} -1} ^{n-1}\sum_{p=j+1}^{\infty}  \frac{ \pi_{p}}{\lambda_{j}\pi_{j}} \ w_j.\]
The second assertion in Lemma \ref{controle} yields the following estimates:
\begin{equation}
\label{bornes}
\|W^{{\scriptscriptstyle 0}}\|_{\ell^\infty(\{n_*{\scriptscriptstyle (K)} -1,n_*{\scriptscriptstyle (K)} ,\ldots\})}\leq C K\ \ ;\ \ 
\|A\|_{\ell^\infty(\{n_*{\scriptscriptstyle (K)} -1,n_*{\scriptscriptstyle (K)} ,\ldots\})}\le  CK.
\end{equation}
We denote by $\Omega'$ the complex disk centered at the origin and of radius $\frac{1}{3CK}$.
Thus, if $\grho\in\Omega'$, the operator $\mathrm{Id} - \grho A$ is invertible and 
$w = (\mathrm{Id} - \grho A)^{-1} \,\grho W^{{\scriptscriptstyle 0}}$. 
It follows from \eqref{bornes} that
\[
\|w\|_{\ell^\infty(\{n_*{\scriptscriptstyle (K)} -1,n_*{\scriptscriptstyle (K)} ,\ldots\})}
\leq \frac{|\grho|\,C K}{1 -|\grho|\,C K}.
\]
Therefore, $w$ is bounded in $\ell^{{\scriptscriptstyle \infty}}(\{\ns-1,\ns,\ldots\})$ by $\frac{1}{2}$ and $1+w_{n}>0$ for all $n\ge \ns-1$. 
It also follows that $w$ is analytic in $\Omega'$. Its derivative is
\[
\frac{\mathrm{d}w}{\mathrm{d}\grho} (\grho) = (\mathrm{Id} - \grho A)^{-1} \,
W^{{\scriptscriptstyle 0}} + (\mathrm{Id} - \grho A)^{-2} \,\grho A\,W^{{\scriptscriptstyle 0}}.
\]
Using \eqref{bornes}, we get for every $\grho\in\Omega'$
\[
\left\|\frac{\mathrm{d}w}{\mathrm{d}\grho}(\grho) -
W^{{\scriptscriptstyle 0}}\right\|_{\ell^\infty(\{n_*{\scriptscriptstyle (K)} -1,n_*{\scriptscriptstyle (K)} ,\ldots\})}\le 4(CK)^2\ |\grho|.
\]
The proof of the proposition is complete.
\end{proof}

\subsection{Matching.}

We consider $I=\left[ -{\scriptstyle 1/(3\widetilde{C}K\log K)},{\scriptstyle 1/(3\widetilde{C}K\log K)}\right]$ and  $K$  large enough so that $\widetilde{C}\log K > C$. With this choice for the interval $I$,
Propositions \ref{ledelta} and \ref{lew} apply for any $\rho\in I$.
We will match the solutions obtained in the two previous subsections in the set $\{\ns-1, \ns\}$,    namely
$u_n^{{\scriptscriptstyle 0}} (1+\delta_n(\grho))$ for $n\leq \ns$ and $1+w_n(\grho)$  for $n\geq \ns-1$.
We will prove that there is a unique  $\grho\in I$ such that there exists a nonzero constant 
$b$ such that  for $n=\ns-1$ and $n=\ns$,
\[
u_n^{{\scriptscriptstyle 0}} (1+\delta_n(\grho)) =  b(1+w_n(\grho)).
\]
We have the following proposition.
\begin{proposition}\label{propprop}
\leavevmode\\
Define the function $f$  by
\begin{align*}
 f(\grho) &=u_{n_{*}{\scriptscriptstyle (K)}-1}^{{\scriptscriptstyle 0}} (1+\delta_{n_{*}{\scriptscriptstyle (K)}-1}(\grho))
(1+w_{n_{*}{\scriptscriptstyle (K)}}(\grho))\\
&  \quad - u_{n_{*}{\scriptscriptstyle (K)}}^{{\scriptscriptstyle 0}}(1+\delta_{n_{*}{\scriptscriptstyle (K)}}(\grho))
(1+w_{n_{*}{\scriptscriptstyle (K)}-1}(\grho)) .
\end{align*}
The minimal positive  zero $\tilde \grho_0$ of $f$ satisfies 
\begin{align*}
& \tilde \grho_0= \frac{\left(1-\frac{\mu_1}{\lambda_1}\right)\sqrt{\frac{\mu_1}{\lambda_1}}\ 
\sqrt{K H''(x_*)}\, x_*\tilde\lambda(x_*)}{\sqrt{2\pi }}
\, \mathlarger{e}^{
K\mathlarger{\int}_{\frac{1}{K}}^{\frac{n_{*}{\scriptscriptstyle (K)}}{K}}
\log \frac{\tilde{\mu}(x)}{\tilde{\lambda}(x)}\,\mathrm{d}x
}\\
& \qquad \times
\mathsmaller{\left(1+ \mathcal{O}\left( \frac{(\log K)^3}{\sqrt{K}}\right)\right)}
\end{align*}
where $H$ is defined in \eqref{def:H}.
\end{proposition}

\begin{proof}
We are going to find a symmetric interval centered around $0$ that contains a unique solution of $f(\grho)=0$. 
Define the auxiliary function $g(\grho)=f(\grho)-f(0)$.
One can check, using Propositions \ref{ledelta} and \ref{lew}  and Lemma \ref{lemma-pouic}  that for all
$\grho\in I$ one has
\begin{align}
&\left| \frac{\mathrm{d} g}{\mathrm{d} \grho}(\grho)-D(K)\right|
 \leq 39 \big( |u_{n_{*}{\scriptscriptstyle (K)}}^{{\scriptscriptstyle 0}}-u_{n_{*}{\scriptscriptstyle (K)}-1}^{{\scriptscriptstyle 0}}| \
\widetilde{C} K \log K + C\widetilde{C}^2 |\grho| K^2 (\log K)^2\big)
\label{gdk}
\end{align}
where 
\begin{equation}\label{def-DK}
D(K)=u_{\ns}^{{\scriptscriptstyle 0}}\big(W_{n_{*}{\scriptscriptstyle (K)}}^{{\scriptscriptstyle 0}} + \Delta_{n_{*}{\scriptscriptstyle (K)}-1}^{\!{\scriptscriptstyle 0}}-\Delta_{n_{*}{\scriptscriptstyle (K)}}^{\!{\scriptscriptstyle 0}}\big)>0.
\end{equation}
Let
\begin{equation}\label{def:etaK}
\eta(K)=
\frac{10}{3}
\frac{u_{n_{*}{\scriptscriptstyle (K)}}^{{\scriptscriptstyle 0}}-u_{n_{*}{\scriptscriptstyle (K)}-1}^{{\scriptscriptstyle 0}}}{D(K)}>0.
\end{equation}
For all $K$ large enough we have, using Lemma \ref{lemma-pouic} items 1 and 4,
\[
\eta(K)\leq \mathcal O(\sqrt{K}) e^{-c K}.
\]
Hence 
\begin{align*}
& \inf_{|\rho|<\eta{\scriptscriptstyle (K)}} \frac{\mathrm{d} g}{\mathrm{d} \grho}(\grho)\\
& \geq D(K)-39\widetilde{C} \big(K \log K |u_{n_{*}{\scriptscriptstyle (K)}}^{{\scriptscriptstyle 0}}-
u_{n_{*}{\scriptscriptstyle (K)}-1}^{{\scriptscriptstyle 0}}|
+ C\widetilde{C} \eta(K) K^2 (\log K)^2\big)\\
& \geq \frac{D(K)}{3}
\end{align*}
for all $K$ large enough by Lemma \ref{lemma-pouic}.
Therefore the function $g$ is monotone increasing in the interval $(-\eta(K),\eta(K))$ and, since $g(0)=0$,
we have
\begin{align*}
& \left[ -\frac{10}{9}(u_{n_{*}{\scriptscriptstyle (K)}}^{{\scriptscriptstyle 0}}-
u_{n_{*}{\scriptscriptstyle (K)}-1}^{{\scriptscriptstyle 0}}), 
\frac{10}{9}(u_{n_{*}{\scriptscriptstyle (K)}}^{{\scriptscriptstyle 0}}\!-u_{n_{*}{\scriptscriptstyle (K)}-1}^{{\scriptscriptstyle 0}})\right] \\
&= \left[ -\frac{ \eta(K) D(K)}{3},\frac{ \eta(K) D(K)}{3}\right]\mathlarger{\subset}\, g\left([-\eta(K),\eta(K)] \right).
\end{align*}
Now because
\[
-f(0)=u_{n_{*}{\scriptscriptstyle (K)}-1}^{{\scriptscriptstyle 0}}-u_{n_{*}{\scriptscriptstyle (K)}}^{{\scriptscriptstyle 0}}\in\!
\left[ -{\scriptstyle \frac{10}{9}}(u_{n_{*}{\scriptscriptstyle (K)}}^{{\scriptscriptstyle 0}}\!-u_{n_{*}{\scriptscriptstyle (K)}-1}^{{\scriptscriptstyle 0}}), {\scriptstyle \frac{10}{9}}(u_{n_{*}{\scriptscriptstyle (K)}}^{{\scriptscriptstyle 0}}\!-u_{n_{*}{\scriptscriptstyle (K)}-1}^{{\scriptscriptstyle 0}})\right]
\]
we have
\[
-f(0) \in g\left([-\eta(K),\eta(K)] \right).
\]
This implies that the equation $g(\grho)=-f(0)$ has a unique solution $\tilde \grho_0$ in $[-\eta(K),\eta(K)]$. 
(This is a special instance of a more general result on quantitative estimates in the inverse function theorem derived in \cite{sotomayor}.)\\

It follows from \eqref{gdk} that for all $\grho\in [-\eta(K),\eta(K)]$
\begin{align*}
& \left|f(\grho)-f(0)-D(K)\grho\right|\\
&\leq 39\ \big( |u_{n_{*}{\scriptscriptstyle (K)}}^{{\scriptscriptstyle 0}}-u_{n_{*}{\scriptscriptstyle (K)}-1}^{{\scriptscriptstyle 0}}|\ \widetilde{C} |\grho|K \log K
+ \frac{1}{2}C\widetilde{C}^2 |\grho|^2 K^2 \log^2K\big)
\end{align*}
which implies that
\begin{align*}
& \left| 
\frac{u_{n_{*}{\scriptscriptstyle (K)}}^{{\scriptscriptstyle 0}}-u_{n_{*}{\scriptscriptstyle (K)}-1}^{{\scriptscriptstyle 0}}}{D(K)}
-\tilde \grho_0
\right|\\
& \leq
\frac{39}{D(K)}\
\left( \left|u_{n_{*}{\scriptscriptstyle (K)}}^{{\scriptscriptstyle 0}}-u_{n_{*}{\scriptscriptstyle (K)}-1}^{{\scriptscriptstyle 0}}\right| 
\widetilde{C} \eta(K) K \log K + \frac{1}{2}C\widetilde{C}^2 \eta(K)^2 K^2 \log^2K\right).
\end{align*}
Using \eqref{def:etaK} and statements 1 and 4 in Lemma \ref{lemma-pouic}, the proposition follows.
\end{proof}

We now end the proof of Theorem \ref{thm:rho0-phi}. We define a sequence $\tilde \varphi$ by
\[
\begin{cases}
\tilde \varphi_{n} = u^{{\scriptscriptstyle 0}}_{n}(1 + \delta_{n}(\tilde\grho_{0})) &\text{for} \, n\leq \ns\\
\tilde \varphi_{n} = b(1 + w_{n}(\tilde\grho_{0})) & \text{for} \, n\geq \ns
\end{cases}
\]
where $\delta_{n}(\tilde\grho_{0})$ and $w_{n}(\tilde\grho_{0})$ are defined in Propositions \ref{ledelta} and \ref{lew}, and
\[
b=\frac{
u^{{\scriptscriptstyle 0}}_{n_{*}{\scriptscriptstyle (K)}}(1 + \delta_{n_{*}{\scriptscriptstyle (K)}}(\tilde\grho_{0}))}{1 + w_{n_{*}{\scriptscriptstyle (K)}}(\tilde\grho_{0})}.
\]
It also follows from these propositions that $\tilde \varphi$ is bounded and hence belongs to $\ell^2(\pi)$. In addition, we get for $n\geq 1$ 
\[
\lambda_{n}\tilde \varphi_{n+1}+\mu_{n}\tilde\varphi_{n-1}\un_{\{n\ge2\}}-(\lambda_{n}+\mu_{n})\tilde\varphi_{n}=-\tilde\grho_{0}\tilde \varphi_{n}.
\]
Let us consider the sequence $(\tilde \varphi^{(k)})_{k\geq 1}$ of elements in $\ell^2(\pi)$ defined by
$\tilde \varphi^{(k)}_{n} = \tilde \varphi_{n} \un_{\{n\leq k\}}$. 
Remark that for all $k\geq 1$, $\tilde \varphi^{(k)}\in \mathscr{D}$.
A straightforward computation leads to
\[
(L\tilde \varphi^{(k)})_{n} + \tilde \grho_{0} \tilde \varphi^{(k)}_{n }=
\begin{cases}
0 & \text{for} \ n<k\\
-\lambda_{k} \tilde \varphi_{k+1 } & \text{for} \ n=k\\
\mu_{k+1}\tilde \varphi_{k} + \tilde \grho_{0} \tilde \varphi_{k+1 }  & \text{for} \ n=k+1\\
0 & \text{for} \ n>k+1.
\end{cases}
\]
Using assumptions \eqref{extinction} and \eqref{eq:pimu}, we can easily prove that
\[
\lim_{k\to \infty} \left(\| \tilde \varphi- \tilde \varphi^{(k)}\|^2_{\ell^2(\pi)} +
\| (L+\tilde \grho_{0}) \tilde \varphi^{(k)}\|^2_{\ell^2(\pi)} \right) = 0.
\]
This implies that $\tilde \varphi\in \mathcal{D}$ and $L\tilde \varphi=-\tilde \grho_{0}  \tilde \varphi$.
By Theorem \ref{autoadj}, the eigenvector $\varphi$ is positive.
Hence it cannot be orthogonal in $\ell^2(\pi)$ to $ \tilde \varphi$ that is strictly positive by Propositions \ref{ledelta}
and \ref{lew}. Since $L$ is self-adjoint, this implies that $\grho_{0} = \tilde \grho_{0}$ and $\varphi = \tilde \varphi$. 

By Assumption \eqref{seconde}, it follows that
$K\int_{0}^{{\scriptscriptstyle \frac{1}{K}}} \log \frac{\tilde{\mu}}{\tilde{\lambda}}(x)\mathrm{d}x =  
\log\frac{\tilde{\mu}}{\tilde{\lambda}}({\scriptstyle \frac{1}{K}}) + {\mathcal O}({\scriptstyle \frac{1}{K}})= 
\log\frac{\mu_{1}}{\lambda_{1}}+ \mathcal{O}({\scriptstyle \frac{1}{K}})$. 
Therefore, using Proposition \ref{propprop} we obtain
\begin{align*}
& \grho_0= \frac{\left(\sqrt{\frac{\lambda_1}{\mu_1}}-\sqrt{\frac{\mu_1}{\lambda_1}}\right)\ 
\sqrt{K\, H''(x_*)}\,x_*\tilde \lambda(x_*)}{\sqrt{2\pi }}\ 
\mathlarger{\mathlarger{e}}^{-K\mathlarger{\int}_{0}^{x_*} \log \frac{\tilde{\lambda}(x)}{\tilde{\mu}(x)}\mathrm{d}x}\\
& \qquad \times \left(1+ \mathcal{O}\left(\frac{(\log K)^3}{\sqrt{K}}\right)\right).
\end{align*}

The estimate for $\varphi$ follows from Propositions \ref{ledelta} and \ref{lew}.

\section{Proof of Theorem \ref{rho1}}\label{proof-rho1}

\subsection{A Poincar\'e inequality}

The proof is based on a Poincar\'e inequality for the Dirichlet form defined for $y\in {\cal D}$ by 
\[
\dir(y)=-\langle y,\gen y\rangle_{\pi}\;.
\]
Recall that $\varphi$ is the eigenvector associated to the maximal eigenvalue $-\grho_0$
of $L$ (see Theorem \ref{autoadj}).

\begin{proposition}
\label{poinca}
\leavevmode\\
For every $y\in \mathcal{D}$ such that $\langle \varphi,y\rangle_{\pi}=0$, we have
\begin{equation}\label{bornedir}
\dir(y)-\grho_{0}\|y\|^{2}_{\pi}\geq g \|y\|^{2}_{\pi}
\end{equation}
where
\begin{equation}
\label{gap}
\frac{1}{g}=\inf_{\tilde{n}\geq 1}
\left(\!\sum_{n=1}^{\tilde{n}}
\frac{1}{\pi_n\lambda_{n}\varphi_{n}\varphi_{n+1}}\sum_{q=1}^{n}\pi_{q}\varphi_{q}^{2}
+\!\sum_{n=\tilde{n}+1}^{\infty}\frac{1}{\pi_{n}\lambda_{n}\varphi_{n}\varphi_{n+1}}\;
\sum_{q=n+1}^{\infty}\pi_{q}\varphi_{q}^{2}\;
\right).
\end{equation}
\end{proposition}

\begin{proof}
Take any $y\in\mathscr{D}$. This implies that there exists some integer $N$
such that $y_n=0$ for all $n>N$. We then have
\begin{align*}
\|y\|^{2}_{\pi}
&=\sum_{n=1}^{N}\bar y_{n}\;\pi_{n}\; y_{n}=
\sum_{n=1}^{N}\frac{\bar y_{n}}{\varphi_{n}}\;y_{n}\pi_{n}\varphi_{n}\\
& =\sum_{n=1}^{N}\frac{\bar y_{n}}{\varphi_{n}}\left(\sum_{q=1}^{n}y_{q}\pi_q\varphi_{q}
-\sum_{q=1}^{n-1}y_{q}\pi_q\varphi_{q}\right)
\end{align*}
where by convention $\sum_{1}^{0}=0$. (Recall that $\bar y_{n}$ is the complex conjugate of $y_{n}$.)
Hence, since $y_{N+1}=0$,
\[
\|y\|^{2}_{\pi}=\sum_{n=1}^{N}\left(\frac{\bar y_{n}}{\varphi_{n}}-
\frac{\bar y_{n+1}}{\varphi_{n+1}}\right)\sum_{q=1}^{n}y_{q}\pi_q\varphi_{q}\;.
\]
By Cauchy-Schwarz inequality we get
\begin{align}
\nonumber
\|y\|^{2}_{\pi}
&
=\sum_{n=1}^{N}\sqrt{\pi_n\lambda_{n}\,\varphi_{n}\,\varphi_{n+1}}
\left(\frac{\bar y_{n}}{\varphi_{n}}-\frac{\bar y_{n+1}}{\varphi_{n+1}}\right)\;
\left(\frac{1}{\sqrt{\pi_n\lambda_{n}\,\varphi_{n}\,\varphi_{n+1}}}\;
\sum_{q=1}^{n}y_{q}\pi_q\varphi_{q}\right)\\
\label{CS}
& \le \sqrt{T_1} \; \sqrt{T_2}
\end{align}
where 
\[
T_1:=\sum_{n=1}^{N}\pi_n \lambda_{n}\,\varphi_{n}\,\varphi_{n+1}
\left|\frac{y_{n}}{\varphi_{n}}-
\frac{y_{n+1}}{\varphi_{n+1}}\right|^{2}
\]
and
\[
T_2:=\sum_{n=1}^{N}\frac{1}{\pi_n \lambda_{n}\,\varphi_{n}\,\varphi_{n+1}}\,
\left|\sum_{q=1}^{n}y_{q}\pi_q\varphi_{q}
\right|^{2}.
\]
Using that $y_{N+1}=0$ and \eqref{tricky} we obtain
\begin{align*}
T_1 & =\sum_{n=1}^{N}
\lambda_{n}\,\pi_n\,\frac{\varphi_{n+1}}{\varphi_{n}}|y_{n}|^{2}
+\sum_{n=1}^{N}\lambda_{n}\,\pi_n\,\frac{\varphi_{n}}{\varphi_{n+1}} |y_{n+1}|^{2} 
-\sum_{n=1}^{N}\lambda_{n}\,\pi_n\,\bar y_{n}y_{n+1}\\
& \quad -\sum_{n=1}^{N}\lambda_{n}\,\pi_n\,\bar y_{n+1}y_{n}\\
& =\sum_{n=1}^{N}
\lambda_{n}\,\pi_n\,\frac{\varphi_{n+1}}{\varphi_{n}} |y_{n}|^{2}
+\sum_{n=2}^{N+1}\lambda_{n-1}\,\pi_{n-1}\,\frac{\varphi_{n-1}}{\varphi_{n}} |y_{n}|^{2}
-\sum_{n=1}^{N}\lambda_{n}\,\pi_n\,\bar y_{n}y_{n+1}\\
& \quad -\sum_{n=1}^{N}\lambda_{n}\,\pi_n\,\bar y_{n+1}y_{n}\\
& =\sum_{n=2}^{N}
\left(\lambda_{n}\,\pi_n\,\frac{\varphi_{n+1}}{\varphi_{n}}+
\lambda_{n-1}\,\pi_{n-1}\,\frac{\varphi_{n-1}}{\varphi_{n}}
\right)|y_{n}|^{2}
+\lambda_{1}\,\pi_{1}\,\frac{\varphi_{2}}{\varphi_{1}}|y_{1}|^{2} \\
&\quad -\sum_{n=1}^{N}\lambda_{n}\,\pi_n\,\bar y_{n}y_{n+1}
-\sum_{n=1}^{N}\lambda_{n}\,\pi_n\,\bar y_{n+1}y_{n}\\
& =\sum_{n=2}^{N}
\left(\lambda_{n}\,\pi_n\,\frac{\varphi_{n+1}}{\varphi_{n}}+
\mu_{n}\,\pi_n\,\frac{\varphi_{n-1}}{\varphi_{n}}
\right)|y_{n}|^{2}
+\lambda_{1}\,\pi_{1}\,\frac{\varphi_{2}}{\varphi_{1}}|y_{1}|^{2}\\
& \quad -\sum_{n=1}^{N}\lambda_{n}\,\pi_n\,\bar y_{n}y_{n+1}
-\sum_{n=1}^{N}\lambda_{n}\,\pi_n\,\bar y_{n+1}y_{n}\\
& =
\sum_{n=1}^{N}(\lambda_{n}+\mu_{n})\,\pi_n|y_{n}|^{2}
-\sum_{n=1}^{N}\lambda_{n}\,\pi_n\,\bar y_{n}y_{n+1}
-\sum_{n=1}^{N}\lambda_{n}\,\pi_n\,\bar y_{n+1}y_{n}\\
& \quad -\grho_{0}\sum_{n=1}^{N}\pi_n\,|y_{n}|^{2}
\end{align*}
since for $n\ge 2$
\[
\lambda_{n}\varphi_{n+1}+\mu_{n}\varphi_{n-1}-(\lambda_{n}+\mu_{n})\,
\varphi_{n}=-\grho_{0}\,\varphi_{n}
\]
and
\[
\lambda_{1}\varphi_{2}-(\lambda_{1}+\mu_{1})\varphi_{1}=-\grho_{0}\varphi_{1}.
\]
Note also that (since $y_{N+1}=0$ and using \eqref{tricky})
\[
\sum_{n=1}^{N}\lambda_{n}\,\pi_n\,\bar y_{n+1}y_{n}
=\sum_{p=2}^{N+1}\lambda_{p-1}\,\pi_{p-1}\,\bar y_{p}y_{p-1}
\]
\[
=\sum_{p=2}^{N+1}\mu_{p}\,\pi_{p}\,\bar y_{p}y_{p-1}
=\sum_{p=2}^{N}\mu_{p}\,\pi_{p}\,\bar y_{p}y_{p-1}.
\]
Therefore
\[
T_1 \le \dir(y)-\grho_{0}\|y\|^{2}_{\pi},
\]
and we get from \eqref{CS} and the previous estimate
\begin{equation}\label{dirdir}
\|y\|^{2}_{\pi} \le  \sqrt{\dir(y)-\grho_{0}\|y\|^{2}_{\pi}}
\; \sqrt{T_2}.
\end{equation}
We now derive an upper bound for $T_2$. We now use the assumption that $y$ is such that
$\langle\varphi,y\rangle_{\pi}=0$ on the top of being such that $y_n=0$ for all $n\geq N+1$. 
In other words
\begin{equation}\label{cancel}
\sum_{q=1}^{N}y_{q}\pi_q\varphi_{q}=
\sum_{q=1}^{\infty}y_{q}\pi_q\varphi_{q}=\langle\varphi,y\rangle_{\pi}=0.
\end{equation}
Let $\tilde{n}$ be a fixed integer over which we will optimize later on.
Then we get, using Cauchy-Schwarz inequality,
\begin{align*}
T_2& \leq  
\sum_{n=1}^{\infty} \frac{1}{\pi_n\lambda_{n}\varphi_{n}\varphi_{n+1}} \left|\sum_{q=1}^{n}y_{q}\varphi_{q}\pi_q
\right|^{2}\\
& =\sum_{n=1}^{\tilde{n}} \frac{1}{\pi_n\lambda_{n}\varphi_{n}\varphi_{n+1}}
\left|\sum_{q=1}^{n}y_{q}\varphi_{q}\pi_q\right|^{2} +
\sum_{n=\tilde{n}+1}^{\infty}\frac{1}{\pi_n\lambda_{n}\varphi_{n}\varphi_{n+1}}
\left|\sum_{q=1}^{n}y_{q}\varphi_{q}\pi_q\right|^{2}\\
& =\sum_{n=1}^{\tilde{n}}
\frac{1}{\pi_n\lambda_{n}\varphi_{n}\,\varphi_{n+1}}
\left|\sum_{q=1}^{n}y_{q}\varphi_{q}\pi_q\right|^{2}
+\sum_{n=\tilde{n}+1}^{\infty} \frac{1}{\pi_n\lambda_{n}\varphi_{n}\varphi_{n+1}}
\left|\sum_{q=n+1}^{\infty}y_{q}\varphi_{q}\pi_q\right|^{2}\\
& \le  \sum_{n=1}^{\tilde{n}}
\frac{1}{\pi_n\lambda_{n}\varphi_{n}\varphi_{n+1}} \left(\sum_{q=1}^{n}|y_{q}|^{2}\pi_q\right) \left(\sum_{q=1}^{n}\varphi_{q}^{2}\pi_q\right)\\
& \quad +\sum_{n=\tilde{n}+1}^{\infty} \frac{1}{\pi_n\lambda_{n}\varphi_{n}\,\varphi_{n+1}}
\left(\sum_{q=n+1}^{\infty}|y_{q}|^{2}\pi_q\right) \left(\sum_{q=n+1}^{\infty}\varphi_{q}^{2}\pi_q\right)\\
& \le
\left(\sum_{n=1}^{\tilde{n}} \frac{1}{\pi_n\lambda_{n}\,\varphi_{n}\varphi_{n+1}}
\sum_{q=1}^{n}\varphi_{q}^{2}\pi_q +\!\!\sum_{n=\tilde{n}+1}^{\infty}\! \frac{1}{\pi_n\lambda_{n}\varphi_{n}\varphi_{n+1}}
\sum_{q=n+1}^{\infty}\!\varphi_{q}^{2}\pi_q\right)\!\! \|y\|^2_{\pi}.
\end{align*}
We used \eqref{cancel} for the second equality, that is,  $\sum_{q=1}^{n}y_{q}\varphi_{q}\pi_q= -\sum_{q=n+1}^{\infty}y_{q}\varphi_{q}\pi_q$.
Combining \eqref{dirdir} and the previous bound we thus get that, if $\langle \varphi,y\rangle_{\pi}=0$, 
\[
\|y\|^{2}_{\pi} \le  \sqrt{\dir(y)-\grho_{0}\|y\|^{2}_{\pi}}\, \frac{1}{\sqrt{g}}\, \|y\|_{\pi}
\]
where $g$ has been defined in \eqref{gap}. This implies \eqref{bornedir}
on $\mathcal{D}$ by closure.
\end{proof}

\subsection{Lower estimate for the spectral gap}\label{subsec:lbsg}

\begin{lemma}
\label{trou}
\leavevmode\\
The spectral gap is bounded below by $g$ defined in \eqref{gap}:
\[
\grho_{1}- \grho_{0} \ge g.
\]
\end{lemma}

\begin{proof}
Let us consider an eigenvector $y\in \mathcal{D}$ with eigenvalue $-\grho_{1}$.
Since $L$ is self-adjoint in $\ell^2(\pi)$, we have $\langle \varphi,y\rangle=0$.
Therefore we get from inequality \eqref{bornedir} in Proposition \ref{poinca}
\[
- \grho_{1} \|y\|^2_{\pi} = \langle y , L y\rangle_\pi \leq  -(\grho_{0}+g)\|y\|^{2}_{\pi}
\]
and the result follows.
\end{proof}

From what precedes, the proof of Theorem \ref{rho1} boils down to prove the following proposition.
\begin{proposition}
\label{prop:g}
\leavevmode\\
For all $K\geq 2$,
$g\geq \frac{\Oun}{\log K}$ where $g$ is defined in \eqref{gap}.
\end{proposition}

Before giving the proof of this proposition, we introduce the following technical quantities.
Let
\begin{equation}\label{eq:xss}
x_{**}=\inf\left\{x\in\real+ : \frac{\tilde{\lambda}(x)}{\tilde{\mu}(x)}<\frac12\right\}.
\end{equation}
Observe that $x_{**}<\infty$ because of \eqref{infini}. Also observe that $x_{*} < x_{**}$ by the assumptions made on the functions $\tilde\lambda$ and $\tilde\mu$. We also define
\begin{equation}\label{eq:nss}
n_{**}{\scriptstyle (K)}=\lfloor x_{**}K\rfloor.
\end{equation}
We will also need to introduce an integer $n_{***}{\scriptstyle (K)}$ that is defined as follows.
By the assumptions made on the functions $\tilde{\mu}$ and $\tilde{\lambda}$ (see \eqref{coeff} and \eqref{eq:xstar}), there exists a number $\theta$ such that
\[
\frac{\tilde{\mu}(0)}{\tilde{\lambda}(0)}<\theta<1.
\]
Thus we can define the following real number (that is strictly smaller than $x_*$).
\begin{equation}\label{def:xtriple}
x_{***}=\sup\left\{x: \frac{\tilde{\mu}(x)}{\tilde{\lambda}(x)}\leq \theta\right\}.
\end{equation}
Then we define the integer
\begin{equation}\label{def:ntriple}
n_{***}{\scriptstyle (K)}=\lfloor x_{***}K\rfloor.
\end{equation}
By definition 
\[
n_{***}{\scriptstyle (K)} \leq \ns \leq n_{**}{\scriptstyle (K)} .
\]
We now turn to the proof of Proposition \ref{prop:g}.

\begin{proof}
From Lemma \ref{propuno} and Theorem \ref{thm:rho0-phi}, we have
\[
\sup_{n\in \entiers\!,K} \left|\varphi_{n} \right| <+\infty\quad\text{and}\quad
\sup_{n\in \entiers\!,K} \left| \varphi_{n}^{-1}\right| <+\infty.
\]
Therefore
\begin{equation}
\label{eq:1surg}
\frac{1}{g}\leq \Oun
\left(\sum_{n=1}^{n_{*}{\scriptscriptstyle (K)}}
\frac{1}{\pi_n\lambda_{n}}\; \sum_{q=1}^{n}\pi_{q}\;
+\sum_{n=n_{*}{\scriptscriptstyle (K)}+1}^{\infty}\frac{1}{\pi_{n}\lambda_{n}}\;
\sum_{q=n+1}^{\infty}\pi_{q}
\right).
\end{equation}

We now derive an upper bound for each sum.\\
We first deal with the second sum in \eqref{eq:1surg}. To this end we write 
\[
\sum_{n=n_{*}{\scriptscriptstyle (K)}+1}^{\infty}\frac{1}{\pi_{n}\lambda_{n}}\;
\sum_{q=n+1}^{\infty}\pi_{q} =
S_{1}+S_{2}+S_{3}
\]
where
\begin{align*}
& S_{1} =
\sum_{n=\nss}^{\infty}\sum_{q=n+1}^{\infty} \frac{\pi_{q}}{\lambda_{n}\pi_{n}} \, ,\quad
S_{2}  = \sum_{n=n_{*}{\scriptscriptstyle (K)}+1}^{\nss}\sum_{q=\nss+1}^{\infty} \frac{\pi_{q}}{\lambda_{n}\pi_{n}} \\
 & \text{and}\quad S_{3}  =
\sum_{n=n_{*}{\scriptscriptstyle (K)}+1}^{\nss}\sum_{q=n+1}^{\nss}
\frac{\pi_{q}}{\lambda_{n}\pi_{n}}
\end{align*}
where $n_{**}{\scriptstyle (K)}$ is defined in \eqref{eq:nss}.
Using Young's inequality and Lemma \ref{controle}, we first get
\[
S_{1}  \le  \sum_{n=\nss}^{\infty}\sum_{q=n+1}^\infty \frac{1}{\mu_{q}}
\left(\frac{1}{2}\right)^{q-n-1}\le \Oun \sum_{q=\nss}^{\infty}  \frac{1}{\mu_{q}} \le \Oun.
\]
Next we have
\begin{align*}
S_{2} & =
\sum_{n=n_{*}{\scriptscriptstyle (K)}+1}^{\nss}\sum_{q=\nss+1}^{\infty} \frac{\pi_{q}}{\lambda_{n}\pi_{n}}
\le\sum_{n=n_{*}{\scriptscriptstyle (K)}+1}^{\nss}\sum_{q=\nss+1}^{\infty}   \frac{1}{\mu_{q}} \left(\frac{1}{2}\right)^{q-\nss-1}\\
&  \le \frac{\Oun}{K} \sum_{n=n_{*}{\scriptscriptstyle (K)}+1}^{\nss}\sum_{q=\nss+1}^{\infty} \left(\frac{1}{2}\right)^{q-\nss-1}
\le\Oun.
\end{align*}
We used several facts: $(\mu_q)$ is increasing, $\tilde{\mu}(x)\geq \tilde{\mu}(0)>0$, and the integers $\ns,\nss$ are of order $K$. 
Finally we have, using Lemma \ref{laborne} and the numbers $\Lambda_{n,m}$ defined just before that lemma,
\begin{align*}
S_{3} & =
\sum_{n=n_{*}{\scriptscriptstyle (K)}+1}^{\nss}\sum_{q=n+1}^{\nss}
 \frac{\pi_{q}}{\lambda_{n}\pi_{n}}
= \sum_{n=n_{*}{\scriptscriptstyle (K)}+1}^{\nss}\sum_{q=n+1}^{\nss}
\frac{\Lambda_{q,n+1}}{\mu_{q}} \\
& \le \frac{\Oun}{K} \sum_{n=n_{*}{\scriptscriptstyle (K)}}^{\nss}\sum_{q=n+1}^{\nss}\;
\mathlarger{e}^{-\,K\,\left(H\left(\frac{q}{K}\right)-H\left(\frac{n+1}{K}\right)\right)}.
\end{align*}
For $x_{*}\le s\le x_{**}$ (see \eqref{eq:xss} for the definition of $x_{**}$) we have for some positive constant $\hat{c}$
\[
\log\frac{\tilde\mu(s)}{\tilde\lambda(s)}\le \hat{c}\;(s-x_{*}).
\]
Hence we get
\begin{align*}
S_{3}
&\le \frac{\Oun}{K} \sum_{n=n_{*}{\scriptscriptstyle (K)}+1}^{\nss}\sum_{q=n+1}^{\nss}\;
\mathlarger{e}^{-\frac{\hat{c}}{2K}\left((q-Kx_{*})^{2}-(n+1-Kx_{*})^{2}\right)}\\
& \le \Oun+ \frac{\Oun}{K}
\sum_{n=n_{*}{\scriptscriptstyle (K)}+1}^{\nss}\sum_{q=n+2}^{\nss}
\;\mathlarger{e}^{-\frac{\hat{c}}{2K}(q-n-1)(q+n+1-2Kx_{*})}
\end{align*}
where we have isolated the term $q=n+1$ that gives $\Oun$. We introduce the new
variables $p=q+n+1$ and $r=q-n-1$ to get
\begin{align*}
S_{3}
& \le \Oun+\frac{\Oun}{K}\sum_{r=1}^{\nss-n_{*}{\scriptscriptstyle (K)}-2}\sum_{p=2n_{*}{\scriptscriptstyle (K)}+r+2}^{2\nss-r}
\;\mathlarger{e}^{-\frac{\hat{c}\,r\,(p-2K x_{*})}{2K}}\\
& =\Oun+\frac{\Oun}{K}\sum_{r=1}^{\nss-n_{*}{\scriptscriptstyle (K)}}\frac{1}{1-\mathlarger{e}^{-\frac{\hat{c} r}{2K}}}\\
& \le \Oun+\frac{\Oun}{K}
\sum_{r=1}^{\nss-n_{*}{\scriptscriptstyle (K)}}\left(\frac{K}{r}+1\right)\\
& \le \Oun+\Oun
\sum_{r=1}^{\nss-n_{*}{\scriptscriptstyle (K)}}\frac{1}{r}\le \Oun\;\log K.
\end{align*}
We now turn to the sum running from $1$ to $\ns$ in \eqref{eq:1surg}. We write
\[
\sum_{n=1}^{n_{*}{\scriptscriptstyle (K)}}\frac{1}{\pi_{n}\lambda_{n}}\;
\sum_{q=1}^{n}\pi_{q}
= \sum_{n=1}^{n_{*}{\scriptscriptstyle (K)}} \sum_{q=1}^{n} \frac{\Lambda_{n+1,q}}{\mu_{q}}  
= \hat{S}_{1}+\hat{S}_{2}+\hat{S}_{3}
\]
where 
\begin{align*}
& \hat{S}_{1}  =\sum_{n=1}^{\nsss}\sum_{q=1}^{n}\frac{\pi_{q}}{\lambda_{n}\pi_{n}}\, ,
\quad\hat{S}_{2} =\sum_{n=\nsss}^{n_*{\scriptscriptstyle (K)}}\sum_{q=1}^{\nsss} \frac{\Lambda_{n+1,q} }{\mu_{q}}\\
& \textup{and}\quad
\hat{S}_{3}=\sum_{n=\nsss}^{n_{*}{\scriptscriptstyle (K)}}\sum_{q=\nsss}^{n}
 \frac{\Lambda_{n+1,q}}{\mu_{q}}.
\end{align*}
By using \eqref{coeff} and inverting the order of summations we get
\[
\hat{S}_{1} \le  \sum_{n=1}^{\nsss} \sum_{q=1}^{n}\frac{1}{\mu_{q}}\;\theta^{n-q}
\le  \sum_{n=1}^{\nsss} \sum_{q=1}^{n}\frac{1}{\tilde{\mu}(0)q }\;\theta^{n-q}
\le \Oun\;\log K
\]
where $n_{***}{\scriptstyle (K)}$ is defined in \eqref{def:ntriple}.
We estimate $\hat{S}_2$ as follows. 
\[
\hat{S}_{2}  \le \Oun \sum_{n=\nsss}^{n_*{\scriptscriptstyle (K)}} \sum_{q=1}^{\nsss} \frac{1}{{q}}  \,\theta^{(\nsss-q)} \leq \Oun .
\]
The last estimate follows by splitting the second sum from $1$ to $\nsss/2$ and from $\nsss/2$ to $\nsss-1$.

Finally, we have the estimates
\[
\hat{S}_{3} \le \frac{\Oun}{K}\!
\sum_{n=\nsss}^{n_*{\scriptscriptstyle (K)}}\sum_{q=\nsss}^{n}
\mathlarger{e}^{K\left(H(\frac{n}{K})-H(\frac{q}{K})\right)}.
\]
For $x_{***}\le s\le x_{*}$ we have
\[
\log\frac{\tilde\mu(s)}{\tilde\lambda(s)}\le -\;c_{2}\;(x_{*}-s)
\]
for some constant $c_{2}>0$, hence
\[
\hat{S}_{3}\le \frac{\Oun}{K}
\sum_{n=\nsss}^{n_*{\scriptscriptstyle (K)}}\sum_{q=\nsss}^{n}\;
\mathlarger{e}^{-\frac{c_{2}}{2K}\,
\left((q-Kx_{*})^{2}-(n-Kx_{*})^{2}\right)}
\]
\[
\le \frac{\Oun}{K}
\sum_{n=\nsss}^{n_*{\scriptscriptstyle (K)}}\sum_{q=\nsss}^{n}\;\mathlarger{e}^{-\frac{c_{2}}{2K}\,(n-q)(2Kx_{*}-q-n)}.
\]
We now use the variables $p=q+n$ and $r=q-n$,
\begin{align*}
\hat{S}_{3}
& \le \Oun+\frac{\Oun}{K}\sum_{r=1}^{n_*{\scriptscriptstyle (K)}-\nsss}\sum_{p=2\nsss+r}^{2n_*{\scriptscriptstyle (K)}-r}
\;\mathlarger{e}^{-c_{2}r\frac{(p-2\,K\,x_{*})}{2K}}\\
& = \Oun+\frac{\Oun}{K}\sum_{r=1}^{n_*{\scriptscriptstyle (K)}-\nsss}\frac{1}{1-\mathlarger{e}^{-\frac{c_{2} r}{2K}}}\\
& \le \Oun+\frac{\Oun}{K}\sum_{r=1}^{n_*{\scriptscriptstyle (K)}-\nsss}\left(\frac{K}{r}+1\right)\\
& \le \Oun+\Oun\sum_{r=1}^{n_*{\scriptscriptstyle (K)}-\nsss}\frac{1}{r}\le \Oun\;\log K.
\end{align*}
Gathering all the bounds, we get the desired result.
\end{proof}

\section{Proof of Theorem \ref{thm:dtv}}\label{proof-thm:dtv}

\subsection{Preliminary estimates}

We first derive some useful estimates. Recall that the constant $c$ has been defined in \eqref{lec}.

\begin{proposition}\label{prop:2bornes}
\leavevmode\\
For all $K>1$ we have
\[
\left|\frac{\langle \varphi,\un\rangle_{\pi}}{\|\varphi\|_{\pi}^{2}}-\frac{1}{u_{n_*{\scriptscriptstyle (K)}}^{{\scriptscriptstyle 0}}}\right|
\leq \Oun\, K^{\frac{3}{2}}\, \log K \, \mathlarger{e}^{-cK}.
\]
\end{proposition}

\begin{proof}
Recall that 
\[
V_n=
\begin{cases}
u_n^{{\scriptscriptstyle 0}} & \textup{if}\; n\leq \ns\\
u_{n_*{\scriptscriptstyle (K)}}^{{\scriptscriptstyle 0}} & \textup{if}\; n\geq \ns.
\end{cases}
\]
Assume that $K$ is large enough so that Propositions \ref{ledelta}, \ref{lew} and Lemma \ref{propuno} apply. We obtain
\begin{equation}\label{eq:phiV}
1-\Oun \grho_0(K) K\log K \leq
\frac{\langle \varphi,\un\rangle_{\pi}}{\|\varphi\|_{\pi}^{2}}
\frac{\|V\|_{\pi}^{2}}{\langle V,\un\rangle_{\pi}}
\leq 1+\Oun \grho_0(K) K\log K.
\end{equation}
Observe that $\|\un\|_\pi^2=\sum_{j=1}^\infty \pi_j$ and
\[
\langle V,\un\rangle_{\pi}-u_{n_*{\scriptscriptstyle (K)}}^{{\scriptscriptstyle 0}} \|\un\|_\pi^2
=\sum_{j=1}^{n_*{\scriptscriptstyle (K)}-1} (u_j^{{\scriptscriptstyle 0}}-u_{n_*{\scriptscriptstyle (K)}}^{{\scriptscriptstyle 0}})\pi_j.
\]
Now using \eqref{def:uno} we get for all $j\leq \ns-1$
\[
u_{n_*{\scriptscriptstyle (K)}}^{{\scriptscriptstyle 0}}-u_j^{{\scriptscriptstyle 0}} = \sum_{\ell=j}^{n_*{\scriptscriptstyle (K)}-1} \frac{1}{\lambda_\ell\pi_\ell}
\]
Hence
\[
\langle V,\un\rangle_{\pi}-u_{n_*{\scriptscriptstyle (K)}}^{{\scriptscriptstyle 0}} \|\un\|_\pi^2=
\sum_{j=1}^{n_*{\scriptscriptstyle (K)}-1}\ \sum_{\ell=j}^{n_*{\scriptscriptstyle (K)}-1} \frac{\Lambda_{\ell+1,j}}{\mu_j}.
\]
We split this sum into three sums, $s_1$, $s_2$ and $s_3$, that we define and estimate as follows.
We have
\[
s_1=\sum_{j=1}^{ \nsss-1}\ \sum_{\ell=j}^{ \nsss-1} \frac{\Lambda_{\ell+1,j}}{\mu_j} \leq \Oun \log K
\]
since in this range $\Lambda_{\ell+1,j}\leq \theta^{\ell-j+1}$ and $\mu_j \geq j \tilde{\mu}(x_*)$.
Next we have
\begin{align*}
s_2 & =\sum_{j=1}^{ \nsss-1}\ \sum_{\ell=\nsss}^{n_*{\scriptscriptstyle (K)}-1} \frac{\Lambda_{\ell+1,j}}{\mu_j} \\
& \leq \Oun \sum_{j=1}^{ \nsss-1}\frac{1}{j} \sum_{\ell=\nsss}^{n_*{\scriptscriptstyle (K)}-1} \Lambda_{\ell+1,\nsss+1}\Lambda_{\nsss,j}.
\end{align*}
We use the fact that $\Lambda_{\ell+1,\nsss+1}\leq 1$ and $\Lambda_{\nsss,j}\leq \theta^{\nsss-j}$ to get
\begin{align*}
s_2 & \leq \Oun \sum_{j=1}^{ \nsss-1}\frac{1}{j} \sum_{\ell=\nsss}^{ n_*{\scriptscriptstyle (K)}-1} \theta^{\nsss-j}\\
& \leq \Oun K \sum_{j=1}^{ \nsss-1}\frac{1}{j}\, \theta^{\nsss-j}\leq \Oun 
\end{align*}
that can be seen by estimating the sums from $1$ to $\nsss/2$ and from $\nsss/2$ to $\nsss-1$.
Finally
\begin{align*}
s_3 & =\sum_{j=\nsss}^{n_*{\scriptscriptstyle (K)}-1}\ \sum_{\ell=j}^{n_*{\scriptscriptstyle (K)}-1} \frac{\Lambda_{\ell+1,j}}{\mu_j} \\
& \leq \Oun\frac{1}{K} \sum_{j=\nsss}^{n_*{\scriptscriptstyle (K)}-1}\ \sum_{\ell=j}^{n_*{\scriptscriptstyle (K)}-1}  \Lambda_{\ell+1,j}
\leq \Oun \log K
\end{align*}
where we first interchange the summations and then follow a very similar argument as in the estimate of $S_3$ in the proof of Proposition \ref{prop:g}.
Therefore we obtain
\begin{equation}\label{coco1}
|\langle V,\un\rangle_{\pi}-u_{n_*{\scriptscriptstyle (K)}}^{{\scriptscriptstyle 0}} \|\un\|_\pi^2| \leq \Oun \log K.
\end{equation}
Now observe that
\begin{align*}
\|V\|_\pi^2 - ({u_{n_*{\scriptscriptstyle (K)}}^{{\scriptscriptstyle 0}}})^2 \| \un\|_\pi^2 & =
\sum_{j=1}^{n_*{\scriptscriptstyle (K)}-1} (({u_j^{{\scriptscriptstyle 0}}})^2-({u_{n_*{\scriptscriptstyle (K)}}^{{\scriptscriptstyle 0}}})^2)\,\pi_j\\
& = \sum_{j=1}^{n_*{\scriptscriptstyle (K)}-1}(u_j^{{\scriptscriptstyle 0}}-u_{n_*{\scriptscriptstyle (K)}}^{{\scriptscriptstyle 0}})(u_j^{{\scriptscriptstyle 0}}+u_{n_*{\scriptscriptstyle (K)}}^{{\scriptscriptstyle 0}})\pi_j.
\end{align*}
Since $(u_j^{{\scriptscriptstyle 0}})$ is monotone increasing and using Lemma \ref{propuno} we get 
\begin{equation}\label{coco2}
|\|V\|_\pi^2 - ({u_{n_*{\scriptscriptstyle (K)}}^{{\scriptscriptstyle 0}}})^2 \| \un\|_\pi^2|\leq \Oun 
\sum_{j=1}^{n_*{\scriptscriptstyle (K)}-1}(u_{n_*{\scriptscriptstyle (K)}}^{{\scriptscriptstyle 0}}-u_j^{{\scriptscriptstyle 0}})\pi_j\leq \Oun \log K
\end{equation}
as we have seen above.\\
Using \eqref{eq:phiV} we have
\begin{align*}
 & \left|\frac{\langle \varphi,\un\rangle_{\pi}}{\|\varphi\|_{\pi}^{2}}-\frac{1}{u_{n_*{\scriptscriptstyle (K)}}^{{\scriptscriptstyle 0}}}\right| \leq
\left|\frac{\langle \varphi,\un\rangle_{\pi}}{\|\varphi\|_{\pi}^{2}}-\frac{\langle V,\un\rangle_{\pi}}{\|V\|_{\pi}^{2}}\right| 
+\left|\frac{\langle V,\un\rangle_{\pi}}{\|V\|_{\pi}^{2}}-
\frac{1}{u_{n_*{\scriptscriptstyle (K)}}^{{\scriptscriptstyle 0}}}\right|\\
& \Oun \grho_0(K) K\log K\, \frac{\langle V,\un\rangle_{\pi}}{\|V\|_{\pi}^{2}}+
\frac{|\,\big\langle V,\un\rangle_{\pi}-u_{n_*{\scriptscriptstyle (K)}}^{{\scriptscriptstyle 0}} \|\un\|_{\pi}^{2}\,\big|}
{\|V\|_{\pi}^{2}}
+
\frac{\big|\,\|V\|_{\pi}^{2}-(u_{n_*{\scriptscriptstyle (K)}}^{{\scriptscriptstyle 0}})^2 \|\un\|_{\pi}^{2}\,\big|}
{u_{n_*{\scriptscriptstyle (K)}}^{{\scriptscriptstyle 0}}\|V\|_{\pi}^{2}}.
\end{align*}
The result follows using \eqref{coco1}, \eqref{coco2}, Lemma \ref{propuno}, and the estimation 
\[
\|V\|_{\pi}^{2}\geq \sum_{n=1}^\infty \pi_n \geq \gamma\, \sqrt{K}\, e^{c K}
\]
where the first inequality follows again from Lemma \ref{propuno} and the definition of $V$, while
the second inequality is the lower bound in statement 5 in Lemma \ref{lemma-pouic}.
\end{proof}

Note that for every $A\in \mathscr{P}(\entiers)$, $\un_{A}\in\ell^{2}(\pi)$. 

\begin{proposition}\label{prop:supA}
\leavevmode\\
There exists $\bar{C}>0$ such that for all $t\geq 0$, for all $K>1$ and for all $n\in\entiers$, we have
\[
\sup_{A\in\mathscr{P}(\entiers)}\left|P_{t}(n,A)-\mathlarger{\mathlarger{e}}^{-\rho_{0}t}
\varphi_{n}\;\frac{\langle \varphi,\un\rangle_{\pi}}{\|\varphi\|_{\pi}^{2}}\, \nu(A)
\right|\le
\frac{\bar{C} K^{\frac{1}{4}} \mathlarger{\mathlarger{e}}^{\frac{c}{2}K} \mathlarger{\mathlarger{e}}^{-\rho_{1}t}}{\sqrt{\pi_{n}}}
\]
where $c$ is defined in \eqref{lec}.
\end{proposition}

\begin{proof}
Let $\mathcal{Q}$ be the spectral projection on the spectral complement of $-\grho_{0}$.
By spectral theory (see {\em e.g.} \cite[Theorem V.2.10, p. 260]{Kato}) we have
\[
\mathlarger{e}^{t L}\un_{A}=
\mathlarger{e}^{-\rho_{0} t}\varphi\;\frac{\langle \varphi,\un_{A}\rangle_{\pi}}{\|\varphi\|_{\pi}^{2}}+
\mathlarger{e}^{t L}\mathcal{Q} \un_{A}\;.
\]
Again by spectral theory  and Cauchy-Schwarz inequality
\[
\left|\langle \mathrm{e}_{n},\mathlarger{e}^{t L}\mathcal{Q}\,\un_{A}\rangle_{\pi}\right|\le 
\mathlarger{e}^{-\rho_{1}t}\;\|\mathrm{e}_{n}\|_{\pi}\; \|\un_{A}\|_{\pi}
\le \mathlarger{e}^{-\rho_{1} t}\,\sqrt{\pi_{n}}\,\sqrt{\sum_{j=1}^\infty \pi_j}.
\]
since $\|\un_A\|_\pi^2\leq \|\un\|_\pi^2 =\sum_{j=1}^\infty \pi_j$. 
The result follows from the definition of $P_t$ (see \eqref{def:Pt}) using statement 5 of Lemma \ref{lemma-pouic}.
\end{proof}

The estimate in Proposition \ref{prop:supA} is not satisfactory for $n$ large since $\pi_n$ tends to
$0$ as $n$ tends to infinity. In fact, we can use the descent from infinity to get an estimate on the error
that is uniform in $n$. 

\begin{proposition}\label{prop:supAmieux}
\leavevmode\\
There exist three strictly positive constants $a,c_1,C'$ such that for all $t\geq 0$, for all $K>1$ and for all $n\in\entiers$,
we have
\[
\sup_{A\in\mathscr{P}(\entiers)}\left|P_{t}(n,A)-\mathlarger{e}^{-\rho_{0}t}
\varphi_{n}\;\frac{\langle \varphi,\un\rangle_{\pi}}{\|\varphi\|_{\pi}^{2}}\, \nu(A)
\right|\le C'\Big( K \mathlarger{e}^{-\frac{a}{4} t} + K^{\frac{3}{4}}\, \mathlarger{e}^{c_1 K} \mathlarger{e}^{-\frac{\rho_1}{2} t}\Big).
\]
\end{proposition}
\begin{proof}
For $q\in\integers$ define $T_q=\inf\{t\geq 0 : X_t^{\scriptscriptstyle{K}}=q\}$.
From the proof of Proposition 2.3 in \cite{BMR} we obtain
\begin{equation}\label{eq:momentexp}
\sup_{n\geq n_{**}{\scriptscriptstyle (K)}}
\esperance_n\left[\mathlarger{e}^{a T_{n_{**}{\scriptscriptstyle (K)}}}\right]\leq \Oun\, K
\end{equation}
where
\[
a=\inf_{K>1}\left(\mathlarger{\mathlarger{\sum}}_{j=n_{**}{\scriptscriptstyle (K)}}^\infty \frac{1}{\lambda_{j}\pi_j}
\mathlarger{\mathlarger{\sum}}_{p=j+1}^\infty \pi_{p}\right)^{-1}.
\]
One can prove that $a>0$ (see Lemma \ref{lemma:a} for a proof).\\
Using Chebyshev inequality we get for all $t>0$
\begin{equation}\label{eq:borneexp}
\sup_{n\geq n_{**}{\scriptscriptstyle (K)}} \proba_n\left(T_{n_{**}{\scriptscriptstyle (K)}}\geq \frac{t}{2}\right)
\leq \Oun \, K \mathlarger{e}^{-\frac{a}{2}t}.
\end{equation}
For every $n\geq n_{**}{\scriptstyle (K)}$, we have
\begin{align*}
P_t(n,A)
&=\esperance_n \left[ \un_A(X_t^{\scriptscriptstyle{K}})\un_{\{T_0>t\}} \right] \\
& =\esperance_n \left[ \un_A(X_t^{\scriptscriptstyle{K}})\un_{\{T_0>t\}} \un_{\{T_{n_{**}{\scriptscriptstyle (K)}}<\frac{t}{2}\}}\right] +
\Oun \; K \mathlarger{e}^{-\frac{a}{2}t}.
\end{align*}
By the strong Markov property we have
\begin{align*}
& \esperance_n \left[ \un_A(X_t^{\scriptscriptstyle{K}})\un_{\{T_0>t\}}
\un_{\{T_{n_{**}{\scriptscriptstyle (K)}}<\frac{t}{2}\}}\right]\\ 
& =
\esperance_n \left[
\esperance_{n_{**}{\scriptscriptstyle (K)}}\!\!\left(\!\un_A\big(X_{t-T_{n_{**}{\scriptscriptstyle (K)}}}^K\big)
\!\un_{\{T_0>t-T_{n_{**}{\scriptscriptstyle (K)}}\}} \!\right)
\!\!\un_{\{T_0>T_{n_{**}{\scriptscriptstyle (K)}}\}} \un_{\{T_{n_{**}{\scriptscriptstyle (K)}}<\frac{t}{2}\}}\right]\!.
\end{align*}
Using Proposition \ref{prop:supA} and Lemma \ref{laborne} we obtain
\begin{align*}
& \esperance_n \!\!\left[ \un_A(X_t^{\scriptscriptstyle{K}})\un_{\{T_0>t\}} \un_{\{T_{n_{**}{\scriptscriptstyle (K)}}<\frac{t}{2}\}}\right]\\
& =
\varphi_{n_{**}({\scriptscriptstyle{K}})} \frac{\langle \varphi,\un\rangle_\pi}{\|\varphi\|_\pi^2}\, \mathlarger{e}^{-\rho_0 t} \nu(A)
\, \esperance_{n} \!\!\left[ \mathlarger{e}^{\rho_0 T_{n_{**}{\scriptscriptstyle (K)}}}\, \un_{\{T_{n_{**}{\scriptscriptstyle (K)}}<\frac{t}{2}\}}\!\right]\!+ \Oun  K^{\frac{3}{4}}\mathlarger{e}^{c_1 K} e^{-\frac{\rho_1}{2} t}
\end{align*}
where $c_1>0$ is a constant independent of $n, t, A$ and $K$. Using Cauchy-Schwarz inequality
we obtain, using \eqref{eq:momentexp} and \eqref{eq:borneexp},
\begin{align*}
\esperance_{n}
\left[ \mathlarger{e}^{\rho_0 T_{n_{**}{\scriptscriptstyle (K)}}}\, \un_{\{T_{n_{**}{\scriptscriptstyle (K)}}\geq \frac{t}{2}\}}\right]
& \leq
\left(\esperance_{n}\left[ \mathlarger{e}^{2\rho_0 T_{n_{**}{\scriptscriptstyle (K)}}}\right]\right)^{\frac{1}{2}}
\left(\esperance_{n}\left[ \un_{\{T_{n_{**}{\scriptscriptstyle (K)}}\geq \frac{t}{2}\}}\right]\right)^{\frac{1}{2}}\\
&\leq \Oun\, K\, \mathlarger{e}^{-\frac{a}{4} t}
\end{align*}
for all $t>0$ and for $K$ large enough so that $2\grho_0\leq a$. Hence
\[
\esperance_{n}
\left[ \mathlarger{e}^{\rho_0 T_{n_{**}{\scriptscriptstyle (K)}}}\, \un_{\{T_{n_{**}{\scriptscriptstyle (K)}}<\frac{t}{2}\}}\right]=\frac{\varphi_n}{\varphi_{n_{**}({\scriptscriptstyle{K}})}}
+ \Oun\, K\, \mathlarger{e}^{-\frac{a}{4} t}
\]
where we used the identity $\esperance_n\left[ e^{\rho_0 T_{n_{**}{\scriptscriptstyle (K)}}}\right]=
\frac{\varphi_n}{\varphi_{\nss}}$ for all $n\geq n_{**}{\scriptstyle (K)}$.
This identity comes from the fact that the process
\[
\left(e^{\rho_{0}(t \wedge T_{n_{**}{\scriptscriptstyle (K)}})} 
\varphi\big(X^{\scriptscriptstyle{K}}_{t\wedge T_{n_{**}{\scriptscriptstyle (K)}}}\big), t\geq 0\right)
\] 
is a martingale (where we write $\varphi(n)$ instead of $\varphi_n$ for the sake of readability).
This relies on the equation $L\varphi = -\rho_{0} \varphi$. The identity then follows from the Martingale Stopping Theorem (see {\em e.g.} \cite{RY}). 
Therefore we obtain
\[
\sup_{A\in\mathscr{P}(\entiers)}\left|P_{t}(n,A)-\mathlarger{e}^{-\rho_{0} t}
\varphi_{n}\frac{\langle \varphi,\un\rangle_{\pi}}{\|\varphi\|_{\pi}^{2}}\nu(A)
\right|
\le \Oun\Big( K \mathlarger{e}^{-\frac{a}{4} t} + K^{\frac{3}{4}}
\mathlarger{e}^{c_1 K} \mathlarger{e}^{-\frac{\rho_1}{2} t}\Big)
\]
for all $n\geq n_{**}{\scriptstyle (K)}$.
The same bound holds for all $n<n_{**}{\scriptstyle (K)}$ using Proposition \ref{prop:supA}.
\end{proof}

\subsection{Proof of Theorem \ref{thm:dtv}}

We first establish inequality \eqref{eq:maindtv}. Observe that for every $B\in\mathscr{P}(\integers)$
\begin{align*}
\proba_n(X_t^{\scriptscriptstyle K} \in B)
&=\proba_n(X_t^{\scriptscriptstyle{K}} \in B\cap \entiers)+ \proba_n(X_t^{\scriptscriptstyle{K}} \in B\cap\{0\})\\
& = \proba_n(X_t^{\scriptscriptstyle{K}} \in B\cap \entiers)+ \gdelta_0(B) \left(1-\proba_n(X_t^{\scriptscriptstyle{K}} \in\entiers)\right).
\end{align*}
Inequality \eqref{eq:maindtv} follows by using twice Proposition \ref{prop:supAmieux}.
This implies the first inequality in the theorem using Proposition \ref{prop:2bornes},  Theorem \ref{thm:rho0-phi} and statement 3 in Proposition \ref{ledelta}.

The second inequality in the theorem is proved as follows.
Let $t_1{\scriptstyle (K)}$ be such that for all $t\geq t_1{\scriptstyle (K)}$
\[
\sup_{n\geq 1} \,\frac{\|\varphi\|_{\pi}^{2}}{\varphi_{n}\langle \varphi,\un\rangle_{\pi}}
\, C'\Big( K\, \mathlarger{e}^{-\frac{a}{4} t} + K^{\frac{3}{4}}\, \mathlarger{e}^{c_1 K} \mathlarger{e}^{-\frac{\rho_1}{2} t}\Big)
\leq \frac12.
\]
We start by considering $t\geq t_1{\scriptstyle (K)}$. We have using Proposition \ref{prop:supAmieux}
\begin{align*}
\left| \frac{P_t(n,A)}{P_t(n,\entiers)}-\nu(A)\right|
&=\left| \frac{P_t(n,A)-\nu(A)P_t(n,\entiers)}{P_t(n,\entiers)}\right|\\
& \leq \frac{2 \|\varphi\|_{\pi}^{2} \, \mathlarger{e}^{\rho_0 t}}{\varphi_n \langle \varphi,\un\rangle_{\pi}}\,
\left| P_t(n,A)-\nu(A)P_t(n,\entiers)\right|.
\end{align*}
The bound follows using again Proposition \ref{prop:supAmieux}, Proposition \ref{prop:2bornes},
Lemma \ref{propuno} (twice), Theorem \ref{thm:rho0-phi} and Propositions \ref{ledelta} and \ref{lew}.
To have the bound for all $t<t_1{\scriptstyle (K)}$, observe that the left-hand side is at most equal to $2$.
The bound follows by eventually taking a larger constant (uniformly in $n$, $K$ and $t$).

\section{Proof of Theorem \ref{thm:ladqs}}\label{proof-thm:ladqs}

Let $K$ be large enough such that $n_{1}=\ns-\sqrt{K}\log K>0$ and
$n_{2}=\ns+\sqrt{K}\log K<n_{**}\scriptstyle{(K)}$.
We have
\[
\frac{\pi_{n}}{\pi_{n_{*}\scriptscriptstyle{(K)}}}=\begin{cases}
\frac{\mu_{n_{*}\scriptscriptstyle{(K)}}}{\mu_{n}}\;\Lambda_{n,n_{*}\scriptscriptstyle{(K)}}^{-1},&\quad\mathrm{for}\quad
n\ge \ns,\\
\frac{\mu_{n}}{\mu_{n_{*}\scriptscriptstyle{(K)}}}\;\Lambda_{n_{*}\scriptscriptstyle{(K)},n},&\quad\mathrm{for}\quad
n\le \ns.
\end{cases}
\]
For $n\le n_{1}$, $\Lambda_{n_{*}\scriptscriptstyle{(K)},n}$ is increasing,
$\mu_{n}\le\Oun K$ and $\mu_{n_{*}\scriptscriptstyle{(K)}}\ge 1$ ($K$ large).
Therefore using Lemma \ref{laborne} we get
\[
\sum_{n=1}^{n_{1}}\frac{\pi_{n}}{\pi_{n_{*}\scriptscriptstyle{(K)}}}\le \Oun K^{2}
\mathlarger{e}^{-c\,(\log K)^{2}}.
\]
Using Lemma \ref{propuno}, Propositions \ref{lew} and \ref{ledelta},
and Theorem \ref{thm:rho0-phi} this implies
\[
\sum_{n=1}^{n_{1}}\frac{\pi_{n}\,\varphi_{n}}{\pi_{n_{*}\scriptscriptstyle{(K)}}\,\varphi_{n_{*}\scriptscriptstyle{(K)}}}
\le \Oun\;K^{2} \mathlarger{e}^{-c\,(\log K)^{2}}.
\]
For  $n_{2}\le n\le n_{**}\scriptstyle{(K)}$, $\Lambda_{n,n_{*}\scriptscriptstyle{(K)}}^{-1}$ is decreasing,
$\mu_{n}\le \Oun K$ and $\mu_{n_{*}\scriptscriptstyle{(K)}}\ge 1$ ($K$ large), therefore using 
Lemma \ref{laborne} we have (since $H''(x_*)>0$)
\[
\sum_{n=n_{2}}^{n_{**}\scriptscriptstyle{(K)}}\frac{\pi_{n}}{\pi_{n_{*}\scriptscriptstyle{(K)}}}
\le \Oun K^{2}\mathlarger{e}^{-c\,(\log K)^{2}}.
\]
For $n\ge n_{**}\scriptstyle{(K)}$ we have
\[
\Lambda_{n,n_{*}\scriptscriptstyle{(K)}}^{-1}\le 
\Lambda_{n_{**}\scriptscriptstyle{(K)},n_{*}\scriptscriptstyle{(K)}}^{-1}\left(\frac{1}{2}\right)^{n-\nss}
\]
hence
\[
\sum_{n=n_{**}\scriptscriptstyle{(K)}}^{\infty}\frac{\pi_{n}}{\pi_{n_{*}\scriptscriptstyle{(K)}}}
\le \Oun K^{2} \mathlarger{e}^{-c\,(\log K)^{2}}.
\]
Using Lemma \ref{propuno}, Propositions \ref{lew} and \ref{ledelta} and Theorem \ref{thm:rho0-phi} this implies
\[
\sum_{n=n_{2}}^{\infty}\frac{\pi_{n}\,\varphi_{n}}{\pi_{n_{*}\scriptscriptstyle{(K)}}\,\varphi_{n_{*}\scriptscriptstyle{(K)}}}
\le\Oun K^{2} \mathlarger{e}^{-c(\log K)^{2}}.
\]
Finally, for $\ns\le n\le n_{2}$, using Lemma \ref{laborne} we have 
\begin{align*}
\frac{\pi_{n}}{\pi_{n_{*}\scriptscriptstyle{(K)}}}
& =\frac{\mu_{n_{*}\scriptscriptstyle{(K)}}}{\mu_{n}}
\sqrt{\frac{\lambda_{n_{*}\scriptscriptstyle{(K)}}}{\mu_{n_{*}\scriptscriptstyle{(K)}}}\frac{\mu_n}{\lambda_n}}\ 
\mathlarger{e}^{-K \left(H\left(\frac{n}{K}\right)-H\left(\frac{n_{*}\scriptscriptstyle{(K)}}{K}\right)\right)- \frac{c(n_{*}\scriptscriptstyle{(K)},n,K)}{K}}\\
& =\mathlarger{e}^{-\frac{(n-n_{*}\scriptscriptstyle{(K)})^{2}}{2K\sigma^{2}}}
\left(1+\frac{\Oun}{K}(n-\ns)+\frac{\Oun}{K^2}(n-\ns)^{3}\right).
\end{align*}
The same estimate holds for $n_{1}\le n\le\ns$.

It is  easy to verify using Lemma \ref{propuno}, Propositions \ref{lew} and \ref{ledelta}, Theorem \ref{thm:rho0-phi} and Lemma \ref{laborne} that for $n_{1}\le n\le n_{2}$
\[
\sup_{n_{1}\le n\le n_{2}}\left|1-\frac{\varphi_{n}}{\varphi_{{n_{*}\scriptscriptstyle{(K)}}}}\right|
\le\frac{\Oun}{K^{2}}.
\]
This implies for $n_{1}\le n\le n_{2}$
\begin{align*}
&\frac{\pi_{n}\,\varphi_{n}}{\pi_{n_{*}\scriptscriptstyle{(K)}}\,\varphi_{n_{*}\scriptscriptstyle{(K)}}}\\
&=\mathlarger{e}^{-\frac{(n-n_{*}\scriptscriptstyle{(K)})^{2}}{2K\sigma^{2}}}
\mathsmaller{\left(1+\frac{\Oun}{K^2}(n-\ns)^{3}\!+\!\frac{\Oun}{K}(n-\ns)+\frac{\Oun}{K^{2}}\right).}
\end{align*}
Therefore, setting
\[
g_{\scriptscriptstyle{K}}(n)=\mathlarger{e}^{-\frac{(n-n_{*}\scriptscriptstyle{(K)})^{2}}{2\,K\,\sigma^{2}}}
\]
we obtain
\[
\sum_{n=n_{1}}^{n_{2}}
\left|\frac{\pi_{n}\,\varphi_{n}}{\pi_{n_{*}\scriptscriptstyle{(K)}}\,
\varphi_{n_{*}\scriptscriptstyle{(K)}}}-g_{\scriptscriptstyle{K}}(n)\right|
\le \Oun.
\]
We also observe that
\[
\sum_{n=1}^{n_{1}}g_{\scriptscriptstyle{K}}(n)+\sum_{n=n_{2}}^{\infty}g_{\scriptscriptstyle{K}}(n)\le 
\Oun \sqrt{K} \mathlarger{e}^{-\tilde c(\log K)^{2}}
\]
for some positive constant $\tilde c$.
Theorem \ref{thm:ladqs} follows after some easy manipulations of the
normalizations.

\section{Appendix: some technical lemmas and estimates}

Let $I=\mathlarger{\int}_{\frac{x_*}{2}}^{+\infty} \frac{\textup{d}x }{x\tilde \mu(x)}$.
Recall that we assume that $I<+\infty$ (see \eqref{cv-int}).

\begin{lemma}\label{controle}
\leavevmode\\
There exists $C\geq 1$ such that for all $K>1$
\[
\sum_{p=n_{*}{\scriptscriptstyle (K)}+1}^\infty \frac{1}{\mu_{p}} \leq I\quad\textup{and}\quad
\sum_{j=n_{*}{\scriptscriptstyle (K)}}^\infty\frac{1}{\lambda_{j}\pi_j} \sum_{p=j+1}^\infty \pi_{p} \ \le \ C\,K.
\]
\end{lemma}

\begin{proof}
Using \eqref{eq:scaling} we get
\[
\sum_{p=n_{*}{\scriptscriptstyle (K)}+1}^\infty \frac{1}{\mu_{p}} =
\frac{1}{K} \sum_{p=n_{*}{\scriptscriptstyle (K)}+1}^\infty \frac{1}{\frac{p}{K}\, \tilde\mu(\frac{p}{K})} 
\leq \int_{\frac{n_{*}{\scriptscriptstyle (K)}}{K}}^\infty \frac{\textup{d}x }{x\, \tilde \mu(x)} \leq I.
\]
This proves the first estimate.
Next, by definition of $\ns$, $n_{**}\scriptstyle {(K)}$ and $x_{**}$ (see Section \ref{sec:hyp} and Subsection \ref{subsec:lbsg}), we have
\begin{align*}
\sum_{j=n_{*}{\scriptscriptstyle (K)}}^\infty \sum_{p=j+1}^\infty \frac{1}{\lambda_{j}} \frac{\pi_{p}}{\pi_{j}}
&= \sum_{j=n_{*}{\scriptscriptstyle (K)}}^{\nss}\sum_{p=j+1}^\infty \frac{1}{\lambda_{j}} \frac{\pi_{p}}{\pi_{j}}+
\sum_{j=\nss+1}^\infty \sum_{p=j+1}^\infty \frac{1}{\lambda_{j}} \frac{\pi_{p}}{\pi_{j}}\\
&\le  \sum_{j=n_{*}{\scriptscriptstyle (K)}}^{\nss}\sum_{p=j+1}^\infty \frac{1}{\mu_{p}}+\sum_{j=\nss+1}^\infty \sum_{p=j+1}^\infty
\left(\frac{1}{2}\right)^{p-j} \frac{1}{\mu_{p}}\\
&\le n_{**}{\scriptstyle (K)}I + \sum_{p=\nss+1}^\infty \frac{1}{\mu_{p}}
\le n_{**}{\scriptstyle (K)} I + \sum_{p=n_{*}{\scriptscriptstyle (K)}+1}^\infty \frac{1}{\mu_{p}}\\
&\leq  (x_{**} \, K + 1)I \leq C K
\end{align*}
where we set $C=(x_{**} + 1)I$ and where we used Young's inequality to get the second inequality.
\end{proof}

\begin{lemma}
\label{lemma:a}
\leavevmode\\
The quantity
\[
a=
\inf_{K>1}\left(\mathlarger{\sum}_{j=\nss}^\infty \frac{1}{\lambda_{j}\pi_j} \mathlarger{\sum}_{p=j+1}^\infty \pi_{p}\right)^{-1}
\]
where $n_{**}\scriptstyle {(K)}$ is defined in \eqref{eq:nss}, is strictly positive.
\end{lemma}

\begin{proof}
The proof follows immediately from the above proof noticing that
\begin{align*}
\sum_{j=\nss+1}^\infty \sum_{p=j+1}^\infty \frac{1}{\lambda_{j}} \frac{\pi_{p}}{\pi_{j}}\leq  \sum_{p=n_{*}{\scriptscriptstyle (K)}+1}^\infty \frac{1}{\mu_{p}}\leq I.
\end{align*}
\end{proof}

Recall that $u^{{\scriptscriptstyle 0}}$ is defined in \eqref{def:uno}.
\begin{lemma}\label{propuno}
\leavevmode\\
There exists a constant $C>0$ such that for all $K$ large enough, and
for all $1\le n\le \ns$
\[
1\le u^{{\scriptscriptstyle 0}}_{n}\leq C.
\]
\end{lemma}
\begin{proof}
We take $K$ large enough such that
\[
1<\lfloor K\,x_{***}\rfloor<\lfloor K\,x_{*}\rfloor-2
\]
where $x_{***}$ is defined in \eqref{def:xtriple}.
Observe that $u_n^{{\scriptscriptstyle 0}}$ is increasing hence for $n\leq \ns$
\begin{align*}
u^{{\scriptscriptstyle 0}}_{n}
& \le u^{{\scriptscriptstyle 0}}_{n_{*}{\scriptscriptstyle (K)}} \\
& \le 1+\mathlarger{\sum}_{j=1}^{n_{***}{\scriptscriptstyle (K)}-1}\frac{1}{\lambda_{j}\,\pi_{j}}
+\mathlarger{\sum}_{j=n_{***}{\scriptscriptstyle (K)}}^{n_{*}{\scriptscriptstyle (K)}}\frac{1}{\lambda_{j}\,\pi_{j}} \\
& \le 
1+ \frac{1}{1-\theta} + \big(\ns -n_{***}{\scriptstyle (K)}+1\big) \theta^{n_{***}{\scriptstyle (K)}}\leq C
\end{align*}
where $C>0$ is independent of $K$.
\end{proof}

For $n> m$ let
\[
\Lambda_{n,m}=\prod_{j=m}^{n-1}\frac{\mu_{j}}{\lambda_{j}}.
\]
By convention we set $\Lambda_{n,n}=1$. We have the following lemma.
\begin{lemma}\label{laborne}
\leavevmode\\
For all $m,n\in\entiers$ such that $n>m$ we have
\[
\Lambda_{n,m}=
\sqrt{\frac{\mu_m}{\lambda_m}\frac{\lambda_n}{\mu_n}}\
\mathlarger{\mathlarger{e}}^{K \left(H\left(\frac{n}{K}\right)-H\left(\frac{m}{K}\right)\right) + \frac{c(m,n,K)}{K}}
\]
where $H$ is defined in  \eqref{def:H} and where $\sup_{m, n, K} |c(m,n,K)| <\infty$.
\end{lemma}

\begin{proof}
By definition \eqref{eq:scaling}
\[
\Lambda_{n,m}=\prod_{j=m}^{n-1}\frac{\tilde{\mu}(\frac{j}{K})}{\tilde{\lambda}(\frac{j}{K})}=
\mathlarger{\mathlarger{e}}^{\sum_{j=m}^{n-1} h\left(\frac{j}{K}\right)}
\]
where $h(s):=\log\big(\tilde{\mu}(s)/\tilde{\lambda}(s)\big)$ ($H'(s)=h(s)$).
Using the trapezoidal rule we get
\begin{align*}
& \log \Lambda_{n,m} \\
&= \frac12 \left(h\left(\frac{n}{K}\right)-h\left(\frac{m}{K}\right)\right) +
K \left(H\left(\frac{n}{K}\right)-H\left(\frac{m}{K}\right)\right) +
\frac{1}{12K^2} \sum_{j=m}^{n-1} h''\left(\frac{\xi_j}{K}\right)
\end{align*}
for some $\xi_j\in [j,j+1]$. Therefore, using \eqref{seconde}, we obtain
\[
\log \Lambda_{n,m} = \frac12 \left(h\left(\frac{n}{K}\right)-h\left(\frac{m}{K}\right)\right)
+ K \left(H\left(\frac{n}{K}\right)-H\left(\frac{m}{K}\right)\right) +  \frac{1}{K}\, c(m,n,K).
\]
The results follows.
\end{proof}

\begin{lemma}
\label{b1}
\leavevmode\\
\[
\textup{For all}\;K>1\;\textup{we have}\;\;u_{n_*{\scriptscriptstyle (K)}}^{{\scriptscriptstyle 0}}=\frac{1}{1-\frac{\mu_1}{\lambda_1}} + \mathcal{O}\left( \frac{1}{K}\right).
\]
\end{lemma}
\begin{proof}
We have
\begin{align*}
u^{{\scriptscriptstyle 0}}_{n_*{\scriptscriptstyle (K)}}
&=1+\mathlarger{\sum}_{j=1}^{n_{*}\scriptscriptstyle{(K)}-1}\frac{1}{\lambda_{j}\pi_{j}}\\
&=1+\mathlarger{\sum}_{j=1}^{\lfloor\sqrt{K}\rfloor}  \prod_{\ell=1}^j \frac{\tilde{\mu}\big(\frac{\ell}{K}\big)}{\tilde{\lambda}\big(\frac{\ell}{K}\big)} + \mathlarger{\sum}_{j=\lfloor\sqrt{K}\rfloor+1}^{n_{*}\scriptscriptstyle{(K)}-1}  \prod_{\ell=1}^j \frac{\tilde{\mu}\big(\frac{\ell}{K}\big)}{\tilde{\lambda}\big(\frac{\ell}{K}\big)}.
\end{align*}
The first sum (plus $1$) is equal to
\begin{align*}
& 1+\sum_{j=1}^{\sqrt{K}}\left(\frac{\tilde{\mu}(0)}{\tilde{\lambda}(0)}\right)^j \mathlarger{e}^{\Oun \frac{j^2}{K}}\\
&= 1+\sum_{j=1}^{\sqrt{K}}\left(\frac{\tilde{\mu}(0)}{\tilde{\lambda}(0)}\right)^j+
\mathcal{O}\left(\frac{1}{K}\sum_{j=1}^{\sqrt{K}} j^2  \left(\frac{\tilde{\mu}(0)}{\tilde{\lambda}(0)}\right)^j\right)\\
&= \frac{1}{1-\frac{\mu_1}{\lambda_1}} + \mathcal{O}\left( \frac{1}{K}\right).
\end{align*}
The second sum is bounded similarly and we get
\[
\Oun K
\left(\frac{\tilde{\mu}(0)}{\tilde{\lambda}(0)}\right)^{\sqrt{{\scriptscriptstyle K}}}\leq\frac{\Oun}{K}.
\]
The lemma is proved.
\end{proof}

The next lemma is about estimating various quantities: $u_{n_*{\scriptscriptstyle (K)}}^{{\scriptscriptstyle 0}}-u_{n_*{\scriptscriptstyle (K)}-1}^{{\scriptscriptstyle 0}}$ (where $u_n^{{\scriptscriptstyle 0}}$ is defined in \eqref{def:uno}), $W_{n_*{\scriptscriptstyle (K)}}^{{\scriptscriptstyle 0}}$ (see \eqref{def-Wn0} for the definition), $\Delta_{n_*{\scriptscriptstyle (K)}}^{\!{\scriptscriptstyle 0}}-\Delta_{n_*{\scriptscriptstyle (K)}-1}^{\!{\scriptscriptstyle 0}}$ (where $\Delta^{\!{\scriptscriptstyle 0}}_n$ is defined in \eqref{def:Deltan0}) and
$D(K)$ (that is defined in \eqref{def-DK}).

\begin{lemma}\label{lemma-pouic}
\leavevmode\\
For all $K>1$ we have the following estimates.
\begin{enumerate}
\item
$u_{n_*{\scriptscriptstyle (K)}}^{{\scriptscriptstyle 0}}-u_{n_*{\scriptscriptstyle (K)}-1}^{{\scriptscriptstyle 0}} = \sqrt{\frac{\mu_1}{\lambda_1}}\ 
\mathlarger{e}^{K\mathlarger{\int}_{\frac{1}{K}}^{\frac{n_*{\scriptscriptstyle (K)}}{K}} \log \frac{\tilde{\mu}(x)}{\tilde{\lambda}(x)}\, \mathrm{d}x}\ \left( 1+ \frac{\Oun}{K}\right);$
\item  
$W_{n_*{\scriptscriptstyle (K)}}^{{\scriptscriptstyle 0}}=
\frac{\sqrt{2\pi}}{2x_* \tilde\lambda(x_*)\sqrt{K H''(x_*)}}\left(1+\frac{(\log K)^3}{\sqrt{K}}\right);$
\item
$\Delta_{n_*{\scriptscriptstyle (K)}}^{\!{\scriptscriptstyle 0}}-\Delta_{n_*{\scriptscriptstyle (K)}-1}^{\!{\scriptscriptstyle 0}}
= -\frac{\sqrt{2\pi}}{2x_* \tilde\lambda(x_*)\sqrt{K H''(x_*)}}  
\left(1+\frac{(\log K)^3}{\sqrt{K}}\right);$
\item
$D(K) = \mathlarger{ \frac{1}{1-\frac{\mu_1}{\lambda_1}}}
\frac{\sqrt{2\pi}}{x_* \tilde\lambda(x_*)\sqrt{K H''(x_*)}}  \left(1+\frac{(\log K)^3}{\sqrt{K}}\right);$
\item

There exist a constant $\gamma\in (0,1)$, that is independent of $K$, such that
\[
\gamma \sqrt{K} \mathlarger{e}^{cK}\leq \sum_{j=1}^\infty \pi_j\leq \gamma^{-1} \sqrt K \mathlarger{e}^{cK},
\]
where $c$ is defined in \eqref{lec}.
\end{enumerate}
\end{lemma}

\begin{proof}
The proof of the first statement follows from Lemma \ref{laborne}, namely
\begin{align*}
u_{n_*{\scriptscriptstyle (K)}}^{{\scriptscriptstyle 0}}-u_{n_*{\scriptscriptstyle (K)}-1}^{{\scriptscriptstyle 0}}
& =\frac{1}{\lambda_{n_*{\scriptscriptstyle (K)} -1}\pi_{n_*{\scriptscriptstyle (K)}-1}}
=\Lambda_{n_*{\scriptscriptstyle (K)},1}\\
& = \sqrt{\frac{\mu_1\lambda_{n_*{\scriptscriptstyle (K)}}}{\lambda_1 \mu_{n_*{\scriptscriptstyle (K)}}}}\ 
\mathlarger{e}^{K\left(H\left(\frac{n_*{\scriptscriptstyle (K)}}{K}\right)-
H\left(\frac{1}{K}\right)\right)+\frac{c(1,n_*{\scriptscriptstyle (K)},K)}{K}}\\
&= \sqrt{\frac{\mu_1}{\lambda_1}}\,
\mathlarger{e}^{K\mathlarger{\int}_{\frac{1}{K}}^{\frac{n_*{\scriptscriptstyle (K)}}{K}} 
\log \frac{\tilde{\mu}(x)}{\tilde{\lambda}(x)}\mathrm{d}x}\, \mathsmaller{\left( 1+ \frac{\Oun}{K}\right).}
\end{align*}

We continue by estimating  $W_{n_{*}\scriptscriptstyle{(K)}}^{{\scriptscriptstyle 0}}$. Write
\begin{align*}
W_{n_*{\scriptscriptstyle (K)}}^{{\scriptscriptstyle 0}}
& = \sum_{p=n_*{\scriptscriptstyle (K)}}^\infty \frac{\pi_p}{\lambda_{n_*{\scriptscriptstyle (K)}-1}\pi_{n_{*}\scriptscriptstyle{(K)}-1}} 
= I_{1}+ I_{2}+I_{3}.
\end{align*}
We start by estimating $I_3$. We again make use of Lemma \ref{laborne}.
\begin{align*}
I_{3}&= \sum_{p=\nss}^\infty \frac{\pi_p}{\lambda_{n_{*}\scriptscriptstyle{(K)}-1}
\pi_{n_{*}\scriptscriptstyle{(K)}-1}} =  \sum_{p=\nss}^\infty  \frac{1}{\mu_p}
\prod_{j=n_{*}\scriptscriptstyle{(K)}}^{\nss-1}\frac{\lambda_j}{\mu_j}\prod_{j=\nss}^{p-1}\frac{\lambda_j}{\mu_j}\\
&\leq \frac{\Oun}{K}\sum_{p=\nss}^\infty \left(\frac{1}{2}\right)^{n-\nss)} \leq \frac{\Oun}{K}
\end{align*}
using the monotonicity of $(\mu_{n})_n$ and the definition of $n_{**}{\scriptstyle (K)}$.\\
We now estimate $I_2$.
\begin{align*}
I_{2}&= \sum_{p=n_{*}\scriptscriptstyle{(K)}+ \sqrt{K}\log K+1}^{\nss-1} \frac{\pi_p}{\lambda_{n_{*}\scriptscriptstyle{(K)}-1}\pi_{n_{*}\scriptscriptstyle{(K)}-1}}=
\sum_{p=n_{*}\scriptscriptstyle{(K)}+ \sqrt{K}\log K+1}^{\nss-1}  \frac{1}{\mu_p\,\Lambda_{p,n_{*}\scriptscriptstyle{(K)}}} \\
&= \sum_{p=n_{*}\scriptscriptstyle{(K)}+ \sqrt{K}\log K+1}^{\nss-1} \frac{1}{\mu_p} 
\sqrt{\frac{\mu_p\lambda_{n_{*}\scriptscriptstyle{(K)}}}{\lambda_p \mu_{n_{*}\scriptscriptstyle{(K)}}}}\ 
\mathlarger{e}^{-K\left(H\left(\frac{p}{K}\right)-H\left(\frac{n_*{\scriptscriptstyle (K)}}{K}\right)\right)}
\mathsmaller{\left( 1+ \frac{\Oun}{K}\right)}\\
&\leq  \sum_{p=n_{*}\scriptscriptstyle{(K)}+ \sqrt{K}\log K+1}^{\nss-1} \frac{1}{\sqrt{\lambda_p \mu_p}} \,
\mathlarger{e}^{-K\left(H\left(\frac{n_*{\scriptscriptstyle (K)}}{K}+ \frac{\log K}{\sqrt{K}}\right)-
H\left(\frac{n_*{\scriptscriptstyle (K)}}{K}\right)\right)}\mathsmaller{\left( 1+ \frac{\Oun}{K}\right)}\\
&\leq \frac{\Oun}{K}(n_{**}{\scriptstyle (K)} -\ns)\,
\mathlarger{e}^{-\frac{(\log K)^2 H''(x_*)}{2} + \frac{\Oun (\log K)^3}{\sqrt{K}}}\\
&\leq \Oun\, \mathlarger{e}^{- \frac{(\log K)^2 H''(x_*)}{2}}
\end{align*}
using the monotonicity of $H$ and Taylor's expansion.\\
Finally we estimate $I_1$. We use again Lemma \ref{laborne}.
\begin{align*}
I_{1}&= 
\mathlarger{\sum}_{p=n_*{\scriptscriptstyle (K)}}^{n_*{\scriptscriptstyle (K)}+ \sqrt{K}\log K} \frac{\pi_p}{\lambda_{n_*{\scriptscriptstyle (K)}-1}\pi_{n_*{\scriptscriptstyle (K)}-1}}\\
&= \mathlarger{\sum}_{p=n_*{\scriptscriptstyle (K)}}^{n_*{\scriptscriptstyle (K)}+ \sqrt{K}\log K} \frac{1}{\mu_p} 
\sqrt{\frac{\mu_p\lambda_{n_*{\scriptscriptstyle (K)}}}{\lambda_p \mu_{n_*{\scriptscriptstyle (K)}}}}\, 
\mathlarger{e}^{-K\left(H\left(\frac{p}{K}\right)-H\left(\frac{n_*{\scriptscriptstyle (K)}}{K}\right)\right)}
\mathsmaller{\left( 1+ \frac{\Oun}{K}\right)}\\
&=\mathlarger{\sum}_{p=n_*{\scriptscriptstyle (K)}}^{n_*{\scriptscriptstyle (K)}+ \sqrt{K}\log K} \frac{1}{\sqrt{\lambda_p \mu_p}}\,
\mathlarger{e}^{-\frac{ H''(x_*) (p-n_*{\scriptscriptstyle (K)})^2}{2 K}}
\mathsmaller{\left( 1+ \frac{\Oun (\log K)^3}{\sqrt{K}}\right)}\\
&= \int_{x_*}^{x_*+ \frac{\log K}{\sqrt{K}}}  \frac{1}{x\sqrt{\tilde\lambda(x)\tilde\mu(x)}} \,
\mathlarger{e}^{-\frac{K H''(x_*)(x-x_*)^2}{2}} \mathrm{d} x \,
\mathsmaller{\left(1+\frac{(\log K)^3}{\sqrt{K}}\right)}\\
&=\frac{1}{2} \frac{\sqrt{2\pi}}{\sqrt{K H''(x_*)}}\frac{1}{x_*\tilde\lambda(x_*) } 
\mathsmaller{\left(1+\frac{(\log K)^3}{\sqrt{K}}\right).}
\end{align*}
The estimation of 
\[
\Delta_{n_*{\scriptscriptstyle (K)}}^{\!{\scriptscriptstyle 0}}-\Delta_{n_*{\scriptscriptstyle (K)}-1}^{\!{\scriptscriptstyle 0}}
=-\sum_{p=1}^{n_*{\scriptscriptstyle (K)}-1}
\frac{(u_p^{{\scriptscriptstyle 0}})^2 \pi_p}{\lambda_{n_*{\scriptscriptstyle (K)}-1}\pi_{n_*{\scriptscriptstyle (K)}-1} 
u_{n_*{\scriptscriptstyle (K)}-1}^{{\scriptscriptstyle 0}} u_{n_*{\scriptscriptstyle (K)}}^{{\scriptscriptstyle 0}}}.
\]
is done similarly by decomposing the sum into three sums with the same ranges as before. \\
The estimation for $D(K)$ follows immediately from the above estimates and Lemma \ref{b1}.\\
Finally, the upper bound in statement 5 is obtained as follows.
We have
\[
\sum_{n=1}^\infty \pi_n = \sum_{n=1}^\infty \frac{1}{\mu_n \Lambda_{n,1}}
\]
where $\Lambda_{n,1}=\prod_{j=1}^{n-1} \frac{\mu_j}{\lambda_j}$.
Using Lemma \ref{laborne} we get
\begin{align*}
\sum_{n=1}^\infty \pi_n &=
\sum_{n=1}^{\nss} \pi_n+ \sum_{n=\nss+1}^\infty \pi_n\\
&\leq \Oun \sum_{n=1}^{\nss} \frac{1}{\sqrt{\lambda_n\mu_n}} \, 
e^{-K\left(H\left(\frac{n}{K}\right)-H\left(\frac{1}{K}\right)\right)}+ \sum_{n=\nss+1}^\infty \pi_n.
\end{align*}
The second sum is estimated by using the fact that $\lambda_j/\mu_j<1/2$ for $j\geq n_{**}{\scriptstyle (K)}$.
The first sum is split into a sum from $1$ to $n_{***}{\scriptstyle (K)}$ and a sum from $n_{***}{\scriptstyle (K)}+1$
to $n_{**}{\scriptstyle (K)}$.
In both cases, we use Lemma \ref{laborne} and the steepest descent method for the 
sum from $n_{***}{\scriptstyle (K)}+1$ to $n_{**}{\scriptstyle (K)}$. The lower bound in statement 5 
is obtained using
\[
\sum_{n=1}^\infty \pi_n \geq \sum_{n=n_{***}{\scriptstyle (K)}}^{\nss} \pi_n
\]
and the steepest descent method as before. This finishes the proof of the lemma.
\end{proof}

Consider the linear equations
\begin{equation}\label{eqabh}
\alpha_{n}w_{n+1}+\beta_{n}w_{n-1}-(\alpha_{n}+\beta_{n})w_{n}=h_{n}
\end{equation}
where $(\alpha_n)_{n\geq1}$, $(\beta_n)_{n\geq1}$ and $(h_n)_{n\geq1}$ are given
sequences of real numbers. The coefficients $\alpha_{n}$ and $\beta_{n}$ are positive.
Define 
\[
\Theta_{p,q}=\prod_{j=q}^{p-1}\frac{\beta_{j}}{\alpha_{j}}\quad\textup{for}\;p>q\quad\textup{and}
\quad \Theta_{q,q}=1.
\]
Note that for $r\ge s\ge  q$
\[
\Theta_{r,s}=\frac{\Theta_{r,q}}{\Theta_{s,q}}.
\]

We have the following lemma.
\begin{lemma}
\label{soleqlin}
\leavevmode\\
The general solution of the homogeneous equation  \eqref{eqabh} when $h_n=0$ for all $n\geq 1$ (homogeneous equation) satisfies the recurrence property
\[
w_{n}=w_{q}+(w_{q+1}-w_{q})\sum_{j=q}^{n-1}\Theta_{j+1,q+1}, \; \forall n\ge q.
\]
In the general case, the solution of \eqref{eqabh} is 
\[
w_{n}=w_{q}+(w_{q+1}-w_{q})\sum_{j=q}^{n-1}\Theta_{j+1,q+1}+
\sum_{j=q}^{n-1} \sum_{p=q+1}^{j}\frac{h_{p}}{\alpha_{p}}\;\Theta_{j+1,p+1}.
\]
In case of convergence of $\sum_{p=q}^{\infty}\frac{h_{p}}{\alpha_{p}\;\Theta_{p+1,q}}$, this can be rewritten as 
\[
w_{n}=w_{q}+\tilde{A}_{q}\sum_{j=q}^{n-1}\Theta_{j+1,q}-
\sum_{j=q}^{n-1} \sum_{p=j+1}^{\infty}\frac{h_{p}}{\alpha_{p}\;\Theta_{p+1,j+1}}, \; \forall n\ge q
\]
for some constant $\tilde{A}_{q}$. \textup{(}We use the convention $\sum_{q}^{q-1}=0$.\textup{)}
\end{lemma}

\begin{proof}
For $n\ge q$ we define $A_{n+1}$ by 
\[
w_{n+1}-w_{n}=A_{n+1}\Theta_{n+1,q}.
\]
Then
\[
\alpha_{n}A_{n+1}\Theta_{n+1,q}-\beta_{n} A_{n}\Theta_{n,q}
=\alpha_{n}\Theta_{n+1,q}\big(A_{n+1}-A_{n}\big)=h_{n}
\]
{\em i.e.}
\[
A_{n+1}-A_{n}=\frac{h_{n}}{\alpha_{n}\Theta_{n+1,q}}
\]
and for all $n\ge q+1$
\[
A_{n}=A_{q}+\sum_{j=q}^{n-1}\frac{h_{j}}{\alpha_{j}\Theta_{j+1,q}}\quad\textup{with}\;
\sum_{q}^{q-1}=0
\]
where 
\[
A_q=\frac{w_{q+1}-w_q}{\Theta_{q+1,q}}-\frac{h_q}{\alpha_q \Theta_{q+1,q} }.
\]
Then for all $n\geq q$
\[
w_{n+1}-w_{n}=A_{q}\Theta_{n+1,q}+\Theta_{n+1,q}\;
\sum_{j=q}^{n}\frac{h_{j}}{\alpha_{j}\Theta_{j+1,q}}
=A_{q}\Theta_{n+1,q}+
\sum_{j=q}^{n}\frac{h_{j}}{\alpha_{j}}\;\Theta_{n+1,j+1}.
\]
Hence
\[
w_{n}=w_{q}+A_{q}\sum_{j=q}^{n-1}\Theta_{j+1,q}+
\sum_{j=q}^{n-1} \sum_{p=q}^{j}\frac{h_{p}}{\alpha_{p}}\;\Theta_{j+1,p+1}.
\]
This implies the first two statements of the lemma.
In case of convergence this can be rewritten as 
\[
w_{n}=w_{q}+\tilde{A}_{q}\sum_{j=q}^{n-1}\Theta_{j+1,q}-
\sum_{j=q}^{n-1} \sum_{p=j+1}^{\infty}\frac{h_{p}}{\alpha_{p}\;\Theta_{p+1,j+1}}
\]
for some constant $\tilde{A}_{q}$. Indeed, since $j\geq p\geq q$, we have $\Theta_{j+1,p+1} = \frac{\Theta_{j+1,q}}{\Theta_{p+1,q}}$. Thus
\begin{align*}
\sum_{j=q}^{n-1} \sum_{p=q}^{j}\frac{h_{p}}{\alpha_{p}}\;\Theta_{j+1,p+1}
&= \sum_{j=q}^{n-1} \Theta_{j+1,q} \sum_{p=q}^{j}\frac{h_{p}}{\alpha_{p}\;\Theta_{p+1,q}}\\
&=  \sum_{j=q}^{n-1}\Theta_{j+1,q}
\left(\sum_{p=q}^{\infty}\frac{h_{p}}{\alpha_{p}\;\Theta_{p+1,q}}
-\sum_{p=j+1}^{\infty}\frac{h_{p}}{\alpha_{p}\;\Theta_{p+1,q}}\right)
\end{align*}
which implies the last statement of the lemma.
\end{proof}

\bigskip \noindent {\bf Acknowledgments}.
The third author benefited from the support of the ``Chaire Mod\'elisation Math\'ematique et Biodiversit\'e''
funded by Veolia Environnement, the Ecole polytechnique and the Mus\'eum national d'Histoire naturelle. The authors thank the referees for their careful reading and comments.
 

\end{document}